\documentclass[preprint,12pt]{elsarticle}

\usepackage{graphicx}
\usepackage{etoolbox}
\makeatletter
\patchcmd{\ps@pprintTitle}{\footnotesize\itshape
       Preprint submitted to \ifx\@journal\@empty Elsevier
       \else\@journal\fi\hfill\today}{\relax}{}{}
\makeatother
\usepackage{bm,amssymb,amsthm,amsmath,mathrsfs}
\usepackage{paralist}
\usepackage{lineno}
\usepackage{fullpage}
\usepackage{mathrsfs}
\usepackage{bm}
\usepackage{natbib}
\usepackage{tabularx}
\usepackage[hidelinks=true]{hyperref}
\hypersetup{
	colorlinks   = true, 
	urlcolor     = blue, 
	linkcolor    = blue, 
	citecolor   = red} 
\usepackage[noprefix]{nomencl}

\renewcommand\nomgroup[1]{%
	\item[\large\bfseries
	\ifstrequal{#1}{Sup}{\it{Superscript}}{%
		\ifstrequal{#1}{Z}{\it{Subscript}}{%
			\ifstrequal{#1}{G}{\it{Greek}}{%
				\ifstrequal{#1}{X}{Other Symbols}}}}]}%
\usepackage{makeidx}         

\usepackage{amssymb}

\usepackage{lineno}

\journal{International Journal of Heat and Mass Transfer \\}

\begin{document}

\begin{frontmatter}



\title{Series Solutions for Orthotropic Diffusion in a Cube}

\author{Brian D. Wood \corref{cor1}\corref{cor2}}
\cortext[cor1]{Corresponding author.}
\cortext[cor2]{e-mail address:brian.wood@oregonstate.edu}
\author{Sassan Ostvar}

\address{School of Chemical, Biological, and Environmental Engineering, Oregon State University, Corvallis, OR  97330 USA}

\begin{abstract}
Analytical solutions to heat or diffusion type equations are numerous, but there are rather few explicit solutions for conditions where the thermal conductivity or diffusion tensors are anisotropic.  Such solutions have some use in making predictions for idealization of real systems, but are perhaps most useful for providing benchmark solutions which can be used to validate numerical codes.  In this short paper, we present the transient solution to the diffusion equation in a cube under conditions of orthotropic anisotropy in the effective thermal conductivity or diffusion tensor.  In particular, we consider the physically-relevant case of transport in a cube with no-flux boundaries for several initial conditions including:  (1) a delta function, (2) a truncated Gaussian function, (3) a step function, and (4) a planar function.   The potential relevance for each of these initial conditions in the context of validating numerical codes is discussed.
 \end{abstract}

\begin{keyword}
\sep diffusion \sep thermal transport \sep analytical solutions \sep heat equation \sep diffusion equation


\end{keyword}

\end{frontmatter}

\setcounter{figure}{0}

\section{Introduction}

Diffusive-like processes in anisotropic media has many applications, ranging from heat transport in composite materials \citep{fischer2003magnetically} to mass and momentum transport in tissues \citep{khaled2003role}.  As a recent topical example, magnetic resonance imaging technologies have been developed to specifically capitalize on the anisotropic structure of neurological tissues; this has lead to a new area of technology known as diffusion tensor imaging (DTI) that has allowed the imaging of neurological tissues at unprecedented levels of detail  \citep{basser1995inferring,basser2011microstructural,basser1994mr,basser1994estimation,westin2002processing,westin1999image}.

Despite the growing interest in diffusive-like processes in anisotropic media, there are very few explicit closed-form solutions available, and fewer that are well-suited for the purposes of testing numerical codes.   The frequently referenced compendium of \citet[][\S1.17-1.20]{carslaw_jaeger_1959} discusses the problem of anisotropy at length, but the only explicit solution presented for a rectangular parallelepiped with constant temperature initial conditions and zero temperature boundary conditions \citep[][\S14.5]{carslaw_jaeger_1959}.  Even then, the solution is left in integral form, and the  coordinate transformation required to achieve a solution is left to the reader to compute.  

In this work, we present explicit series solutions for anisotropic heat transport or diffusion in a spatially homogeneous cube for three different initial condition cases.  The initial conditions have been selected with the particular purpose of being conditions that are interesting from the perspective of validating numerical codes.  Because heat transport and diffusion in a homogeneous anisotropic medium are mathematically identical, the solutions developed in this paper apply to either case.  To be specific, note that the diffusion and heat equations are given by 

\begin{align}
\frac{\partial c_{A\gamma}}{\partial t} &=\nabla\cdot\left(  {\textbf{\sffamily \bfseries D}}_{\gamma}\cdot \nabla c_{A\gamma} \right) \\
\rho c_p \frac{\partial T_{A\gamma}}{\partial t} &=\nabla\cdot\left(  {\textbf{\sffamily \bfseries K}}_{\gamma}\cdot \nabla T_{A\gamma} \right) 
\end{align}
For the case of homogeneous (but anisotropic) systems, we clearly have a correspondence between the diffusion tensor, ${\textbf{\sffamily \bfseries D}}_{\gamma}$, and the thermal conductivity tensor divided by the product of the density and the heat capacity, $ {\textbf{\sffamily \bfseries K}}_{\gamma}/(\rho c_p)$.  In the remainder of this paper, we will present the solutions for the diffusion equation, with the idea that all of the results are just as applicable to heat transport using the simple correspondence discussed.

In this work, we consider the problem of diffusion of a dilute chemical species in a homogeneous but anisotropic medium, whose diffusion tensor can be expressed by a diagonal tensor with non-equal entries (i.e., an orthotropic tensor) .  Our intent is to provide new solutions that focus on novel initial value functions; to date, most of the closed-form solutions have been boundary driven.  Thus, we have developed solutions for initial conditions that are (1) not constant, and (2) of interest for validation of numerical codes.  The novel features of the material presented are primarily the initial conditions covered, and the explicit, closed form of the solutions.

\section{Background}

For the case of diffusion in isotropic systems, a large number of solutions covering many possible geometries, boundary condition, and initial conditions can be found in the classic texts by \citet{carslaw_jaeger_1959} and  \citet{crank1975mathematics}.  Excellent reviews of more modern papers on this topic can be found in references \citep{de2003unsteady} and \citep{lu2005new}.  For a particular set of boundary and initial conditions, we follow  {\it explicit} solutions as those which contain only sums and products of well-known solutions (including conventional algebraic, elementary transcendental, and conventional transcendental special functions); following \citet{olver2014introduction}, such solutions include infinite series.  Series solutions can be very useful when they converge reasonably rapidly because they provide a method to compute high-accuracy approximations to the solution with a finite number of simple computations.  Solutions that contain one or more unevaluated integrals are terms {\it integral} solutions.  Although these solutions are also useful, it may not be possible to evaluate the integrals that are encountered in terms of algebraic, transcendental, or special functions.  Such solutions, ultimately, would require a numerical quadrature to evaluate.

Although the literature on the anisotropic case is reasonably well developed, this literature has focussed primarily (although not exclusively) on (1) systems that are orthrotropic (i.e., the diffusion tensor is diagonal, but the diagonal entries are not equal), and (2) systems that are driven by boundary conditions rather than initial conditions.    The literature for both steady and transient solutions in both 2- and 3- dimensions is summarized in the following paragraphs.   

In 2-dimensions, explicit closed-form steady-state solutions have been generated for both the orthotropic case \citep{akoz1978thermoelastic,tauchert1974thermal}, and for the fully anisotropic case \citep{ma_2004}.  
For the 2-dimensional transient case, separation of variables has been successful for developing explicit closed-form solutions for orthotropic cases \citep[][Ex.~15-6]{sugano1980transient,wang1985transient,milovsevic2004analytical,hahn2012heat}.  The paper by \citet{wang1986transient} is unique because it is the only result that allows for general initial conditions; all others have been been developed for the zero or (equivalently) constant initial condition.  A single paper by \citet{hsieh2002analytical} reports explicit closed-form results for a fully anisotropic system which are specific to a two-layer infinite slab.

In 3-dimensions, there are explicit closed-form solutions for the steady orthotropic case given by \citet{tungikar1994three}, and \citet[][Ex.~15-3, 15-4]{hahn2012heat}.  Aside from the well-known solution for the infinite domain given by \cite{carslaw_jaeger_1959}, the only other fully transient explicit closed-form solution we were able to find in the literature was the solution of \cite{ootao2007three} for a finite orthotropic functionally graded plate with zero initial conditions.  Although \citet{padovan1975solution} outlines a variant of the separation of variables approach for a fully anisotropic transient solution, no explicit solutions are provided, and the results are restricted to particular geometries.  To date, there appears to be no general approach that yields explicit closed-form solutions for fully anisotropic tensors in which the principle axes are not perpendicular to boundaries of the domain.

Integral solutions (either using Green's functions or Laplace transforms) have also been developed by several researchers for both the orthotropic and the fully anisotropic cases in  2-dimensions \citep{aviles1998exact,delouei2012exact,frankel2010new,clements1973thermal} and 3-dimensions \citep{beck2004solutions,chang1977heat,haji2002temperature,marczak_2011,poon_1979,taheri2014analytical,taheri2013transient,chang_1977}.  These solutions are formally explicit, but they are not in a closed form. Thus, they these have the disadvantage of requiring the computation of integrals before an explicit solution can be produced.

\section{Symmetry and Positivity of the Diffusion and Thermal Conductivity  Tensors}
%
In 1851 George Stokes (of fluid dynamics fame) published a paper arguing that thermal transport was tensorial in nature, and that the associated transport tensor was symmetric.  Since that time, there has been substantial discussion on the requirements for the components of a conductivity or diffusion tensor to be physically meaningful (see extensive discussions on the topic in references \citep[][]{Truesdell:1960aa,coleman1960reciprocal,day1969symmetry,miller_1960,powers2004necessity,Truesdell:1960aa,wei1966irreversible}.)   The diffusion tensor (or, equivalently, the thermal conductivity tensor) is a positive definite tensor if one accepts the axioms of reciprocity laid out by Onsager \citep{onsager1931reciprocal_II,onsager1931reciprocal_I,miller_1960}.    Unlike isotropic tensors, anisotropic tensors have the potential to have zero transport in one or more principle directions, and thus can only be assured to be positive semidefinte (rather than positive definite).   Assuming that the tensor takes the form 

\begin{equation}
{D_{ij,\gamma }} = \left[ {\begin{array}{*{20}{c}}
a&d&e\\
d&b&f\\
e&f&c
\end{array}} \right]
\end{equation}
then positive semidefinteness is assured by checking Sylvester's criterion \citep{gilbert1991positive} that all of the principal minors are all nonnegative, i.e., $a\ge 0, b\ge0, c\ge0, ab-d^2\ge0, bc-f^2\ge0$, $ac-e^2\ge0$, and $a(bc- f^2) - d(dc-ef) + e(df-be)\ge 0$ .  At least one of these quantities must be nonzero to prevent the problem from being the trivial case where no transport occurs.   Note, however, that it is not necessary that all of the principal minors be nonzero.  For example, in a composite medium constructed of impermeable matrix and spanning capillaries aligned in the horizontal plane, the effective diffusion coefficient in the vertical direction is zero, and this would be reflected in the principal minors. 

Every three-dimensional second-order positive semidefinite tensor has the following properties: (1) there exist three nonnegative eigenvalues for the tensor; (2) if the three Eigenvalues are distinct and positive, then there exist three unique eigenvectors defining the principal directions of the tensor; (3) these three principle directions are mutually orthogonal, (4) for the case where one or two of the Eigenvalues is zero, or there are positive but repeated Eigenvalues, then there are an infinite number of such mutually orthogonal coordinate systems, and (5) the diffusion tensor is diagonal (but not necessarily isotropic) when expressed in a mutually orthogonal coordinate system constructed from its Eigenvectors.  Thus, the anisotropic diffusion equation given by Eqs.~(\ref{aniso_1})-(\ref{aniso_3}) can always be expressed in terms of a orthotropic (diagonal) tensor, $ {\textbf{\sffamily \bfseries D}}_{\gamma}$, when put in the proper mutually-orthogonal coordinate system; the tensor takes the form
\begin{equation}
{D_{ij,\gamma }} = \left[ {\begin{array}{*{20}{c}}
D_{xx}&0&0\\
0&\tfrac{D_{xx}}{d^2_{yy}}&0\\
0&0&\tfrac{D_{xx}}{d^2_{zz}}
\end{array}} \right]
\label{ortho}
\end{equation}
Here $D_{xx},~D_{yy}$, and $D_{zz}$ represent the diffusion coefficients in the direction of the $x-$, $y-$, and $z-$directions, respectively. The other two constants represent the anisotropy in the system, and are defined by $d^2_{yy}=D_{xx}/D_{yy}$, and $d^2_{zz}=D_{xx}/D_{zz}$. 

For such conditions, the governing differential equation for anisotropic diffusion is well understood.  The problem in an {\it unbounded domain} in 3-dimensions can be stated by

\begin{align}
&&\frac{\partial c_{A\gamma}}{\partial t} &=\nabla\cdot\left(  {\textbf{\sffamily \bfseries D}}_{\gamma}\cdot \nabla c_{A\gamma} \right) \label{aniso_1}\\
&Maximum~Principle&  c({\bf x},t) & \rightarrow 0 ~as~\|{\bf x}\| \rightarrow \infty \\
&Constraint && \nonumber\\
&I.C.& c({\bf x},0)&=\Phi({\bf x})
\label{aniso_3}
\end{align}
This problem requires that the initial condition have finite total mass (or energy); in the following we will assume for convenience that the initial condition has a total mass of unity. 

The problem diffusive-like transport in an infinite medium was introduced first by \citet{stokes_1851}; full solutions to the problem were introduced by Carslaw and Jaeger, although not in explicit form.   However, following their discussions in \S1.17, \S2.2, and \S10.2, one can, after a bit of algebra, extract the solutions for the whole space of the form

\begin{align}
c_{A\gamma}({\bf x},t) =&  \frac{1}{2 \sqrt{\pi D_{xx} t}} \frac{d_{yy}}{2 \sqrt{\pi D_{xx} t}} \frac{d_{zz}}{2 \sqrt{\pi D_{xx} t}}
\int_{x=-\infty}^{x=\infty}  \int_{y=-\infty}^{y=\infty}  \int_{z=-\infty}^{z=\infty}   
\Phi(x',y',z')
\exp\left(  - \frac{(x-x')^2}{4 D_{xx} t} \right) \nonumber \\
\times & \exp\left(  - \frac{(y-y')^2 d^2_{yy}}{4 D_{xx} t} \right)
\exp\left(  - \frac{(z-z')^2 d^2_{zz}}{4 D_{xx} t} \right)\;
dx'\;dy'\;dz'
\end{align}
For the initial condition $ \Phi(x',y',z')=\delta(x')\delta(y')\delta(z')$ this gives an explicit non-integral solution of the form.

\begin{align}
c_{A\gamma}({\bf x},t) =&  \frac{1}{2 \sqrt{\pi D_{xx} t}} \frac{d_{yy}}{2 \sqrt{\pi D_{xx} t}} \frac{d_{zz}}{2 \sqrt{\pi D_{xx} t}}
\exp\left(  - \frac{x^2}{4 D_{xx} t} \right) \exp\left(  - \frac{y^2 d^2_{yy}}{4 D_{xx} t} \right)
\exp\left(  - \frac{z^2 d^2_{zz}}{4 D_{xx} t} \right)\;
\end{align}
which is the frequently referenced solution for an instantaneous point source presented by \citet[][\S 10.2, page 257]{carslaw_jaeger_1959}.

In bounded domains, the variable transformations discussed above can still put the diffusion tensor in orthotropic form; however, the boundary conditions are not, in general, easy to handle (e.g., separable) in this coordinate system.  Thus, coordinate transformation does not lead to tractable solutions for the diffusion problem except when the principal axes of the diffusion tensor happen to align with the boundaries of the domain (i.e., the problem is orthotropic because the domain of the parallelepiped is constructed so that the sides aligned with the tensor principle axes.)  

Because much of the terminology and early work on this topic came originally from research on the properties of crystals, it is worthwhile to think about an example from a crystal system.  A monoclinic crystal is a rhomboid where two of the three planes defining the system are orthogonal, and the third is not. The natural coordinate system describing the geometry of a monoclinic crystal (the crystallographic coordinate system) is one where two of the directional vectors are perpendicular, and the third makes an angle other than $90^{\circ}$.  In the natural coordinate system, the thermal conductivity or diffusion tensor can be expressed by the anisotropic tensor \citep[][\S 1.17]{carslaw_jaeger_1959}

\begin{equation}
{D_{ij,\gamma }} = \left[ {\begin{array}{*{20}{c}}
a&d&0\\
d&b&0\\
0&0&c
\end{array}} \right]
\label{monclinic}
\end{equation}

Although an orthotropic representation of the thermal conductivity or diffusion tensor in such a crystal does have an orthotropic representation in some coordinate system, that coordinate system is not aligned, in general, with the natural coordinate system of the crystal.  Hence, the natural boundaries of the crystal would be complicated functions of more than one variable (and, hence, unlikely to be easily separable), making it difficult to apply separation of variables to the problem directly.  This is not a merely academic issue; experimental measurements have been made for diffusion in crystals that show exactly this behavior.  A paper by \citet{bendani1993diffusion} describes a set of measurements of the diffusion tensor for a naphthalene crystal, and the result is expressed both as an fully anisotropic tensor in the crystallographic coordinate system (i.e., analogous to Eq.~(\ref{monclinic})), and as an orthotropic tensor (analogous to Eq.~\ref{ortho}) when put in the principle axes of the diffusion process.

\section{Orthotropic Diffusion in a Cubic Domain}

In the remainder of this work, we focus specifically on systems where the coordinate system is organized to yield an is orthotropic diffusion tensor, and the boundary conditions are defined on planes perpendicular to the coordinate axes.   We further assume that the domain is a cube with side $L$; general rectangular parallelepipeds domains can always be transformed to a cubic domain as described in the Appendix.  Although a few solutions do exist in the literature, these solutions have been, to date, either 2-dimensional, or for particularly simple initial and boundary conditions.  In this work, we illustrate how the solution can be obtained for arbitrary (assuming only that the function is an admissible one so that the partial differential equation to be well defined) with no-flux conditions at the boundaries. The no-flux conditions correspond closely to conditions that are of interest both for comparison with experimental or numerical results.  For this case, the basic mass balance equation and boundary conditions take the form

\begin{align}
&&\frac{\partial c_{A\gamma}}{\partial t} &=\nabla\cdot\left(  {\textbf{\sffamily \bfseries D}}_{\gamma}\cdot \nabla c_{A\gamma} \right) \\
&B.C.~1& -{\bf n}_{\gamma} \cdot  {\textbf{\sffamily \bfseries D}}_{\gamma}\cdot \nabla c_{A\gamma}  &= 0,~~ at ~ all ~ bounaries~of~cube \\
&I.C.& c({\bf x},0)&=\Phi({\bf x})
\end{align}
where ${\bf L}=(L,L,L)$.  Note that in the following, we will assume that $ {\textbf{\sffamily \bfseries D}}_{\gamma}$ is a diagonal (but anisotropic) tensor.
We will attempt to find a solution using separation of variables.  Assume $c_{A\gamma}(x,y,z,t) =T(t)U(x,y,z)$, where $U(x,y,z)=X(x)Y(y)Z(z)$.  Substituting this proposed solution yields
\begin{align}
\frac{\partial T}{\partial t} U = T  {\textbf{\sffamily \bfseries D}}_{\gamma}: \nabla \nabla U
\end{align}
where here $\nabla\nabla U = \frac{\partial^2 U}{\partial x_i \partial x_j}$.   Note that because $ {\textbf{\sffamily \bfseries D}}_{\gamma}$ is diagonal, we have  ${\textbf{\sffamily \bfseries D}}_{\gamma} = 0~for ~ i \ne j$.   Thus, we have
\begin{equation}
 {\textbf{\sffamily \bfseries D}}_{\gamma}:\nabla\nabla U = D_{i j} \frac{\partial^2 U}{{\partial x_j}{\partial x_i}}= D_{xx}\frac{\partial^2 U}{\partial x^2}+D_{yy}\frac{\partial^2 U}{\partial y^2}+D_{zz}\frac{\partial^2 U}{\partial z^2}
\end{equation}
where each of $D_{xx}, D_{yy},$ and $D_{zz}$ are (potentially) unique.

\subsection{Separation of Variables}
%
For this problem, the process of separation of variables can be used efficiently to develop fully transient, 3-dimensional solutions.  Although the approach is reasonably straight-forward, the separation does require some careful handling of the separation constants involved.  In any event, the separation of this problem does not appear to be well described in the literature.  As a result, the process is described in detail in this section.  

Using the functional relationship $U(x,y,z)=X(x)Y(y)Z(z)$ and substituting, we find

\begin{equation}
XYZ\frac{\partial T}{\partial t} = YZD_{xx}\frac{\partial^2 X}{\partial x^2}T+XZD_{yy}\frac{\partial^2 Y}{\partial y^2}T+XYD_{zz}\frac{\partial^2 Z}{\partial z^2}T
\end{equation}
Dividing both sides of this equation by $XYZT$ yields
\begin{equation}
\frac{1}{T}\frac{\partial T}{\partial t} = \frac{D_{xx}}{X}\frac{\partial^2 X}{\partial x^2}+\frac{D_{yy}}{Y}\frac{\partial^2 Y}{\partial y^2}+\frac{D_{zz}}{Z}\frac{\partial^2 Z}{\partial z^2}
\end{equation}
For later convenience, we normalize this equation by $D_{xx}$.  Defining $d^2_{yy} = D_{xx}/D_{yy}$ and $d^2_{zz} = D_{xx}/D_{zz}$, we have
\begin{equation}
\frac{1}{D_{xx} T}\frac{\partial T}{\partial t} = \frac{1}{X}\frac{\partial^2 X}{\partial x^2}+\frac{1}{d^2_{yy}Y}\frac{\partial^2 Y}{\partial y^2}+\frac{1}{d^2_{zz}Z}\frac{\partial^2 Z}{\partial z^2}
\end{equation}
Noting that for this equation, the right-hand side is a function of only $X,Y,$ and $Z$, and the left-hand side is a function of only $T$, we conclude, as conventional, that both sides must be equal to a constant.  The same argument can be made for each of the terms on the right-hand side independently; each term is independent of all of the others.  Thus, we define the separation constant $\lambda^2$ by
\begin{equation}
\frac{1}{D_{xx} T}\frac{\partial T}{\partial t} = \frac{1}{X}\frac{\partial^2 X}{\partial x^2}+\frac{1}{d^2_{yy}Y}\frac{\partial^2 Y}{\partial y^2}+\frac{1}{d^2_{zz}Z}\frac{\partial^2 Z}{\partial z^2} = -\lambda^2
\end{equation}
we note also that, because each of the terms on the right hand side is independent of all others, that each of these terms is also equal to a constant.  For convenience, we define these three constants by 
\begin{equation}
\lambda^2 = \alpha^2 + \beta^2 + \gamma^2
\end{equation}

This defines a set of four ordinary differential equations as follows

\begin{align}
\frac{d T}{d t} +\lambda^2 D_{xx} T &= 0 \\
\frac{d^2 X}{dx^2} +\alpha^2 X &= 0 \\
\frac{d^2 Y}{dy^2} + d^2_{yy}\beta^2 Y &= 0 \\
\frac{d^2 Z}{dz^2} + d^2_{zz}\gamma^2 Z &= 0 
\end{align}

The solutions to each of these is straightforward.  For the time variable, the solution is a decaying exponential of the form

\begin{equation}
T(t) = T_0 \exp\left( -\lambda^2 D_{xx} t  \right)
\end{equation}

For the remaining equations, the characteristic equation yields roots of the form $\pm \alpha \imath $, $\pm d_{yy} \beta \imath $, and $\pm d_{zz} \gamma \imath $.  Thus, the solutions are trigonometric functions of the form

\begin{align}
X(x) &= A \sin (\alpha x) + B \cos (\alpha x) \\
Y(y) &= G \sin (d_{yy} \beta y) + H \cos (d_{yy} \beta y) \\
Z(Z) &= R \sin (d_{zz} \gamma z) + S \cos (d_{zz} \gamma z) 
\end{align}
The use of the no-flux boundary conditions require that these functions obey the relationship

\begin{align}
X'(x) = \alpha A \cos (\alpha x) - \alpha/L B \sin (\alpha x) \\
Y'(y) = d_{yy} \beta G \cos (d_{yy} \beta y) - d_{yy} \beta H \sin (d_{yy} \beta y) \\
Z'(z) = d_{zz} \gamma R \cos (d_{zz} \gamma z) - d_{zz} \gamma S \sin (d_{zz} \gamma z) 
\end{align}
For each of the three boundaries where the origin is included, these relationships indicate that $A = G = R = 0$.  Thus the solution is of the form

\begin{align}
X(x) & =  B \cos (\alpha x) \\
Y(y) & =  H \cos (d_{yy} \beta y) \\
Z(Z)  & =  S \cos (d_{zz} \gamma z)
\end{align}
The remaining boundary conditions require the following
\begin{align}
\sin (\alpha L) & = 0 \\
\sin (d_{yy} \beta L) & = 0 \\
\sin (d_{zz} \gamma L) & = 0
\end{align}
and these are met, respectively, by Eigenvalues of the form $\alpha = \ell \pi / L$, $d_{yy} \beta = m \pi /L$, and $d_{zz} \gamma = n \pi /L$ for $\ell, m, n = 0, 1, 2 \ldots$  The solution of the set of differential equations, then, is given by the linear combination of all possible Eigenfunctions.  

\begin{equation}
c_{A\gamma}(x,y,z,t) = \sum_{\ell = 0}^{\ell=\infty} \sum_{m = 0}^{n=\infty}\sum_{n = 0}^{m=\infty}   \exp\left( - \lambda^2 D_{xx} t \right)   B_{\ell} H_{m} S_{n} cos(\ell \pi \tfrac{x}{L}) \cos(m \pi \tfrac{y}{L}) \cos(n \pi \tfrac{z}{L})
\end{equation}
Recalling the relationship $\lambda^2= \alpha^2+\beta^2+\gamma^2$, then we also have $\lambda^2 = ( \ell^2  +m^2/d^2_{yy}  + n^2/d^2_{zz} ) \frac{\pi^2}{L^2}$.  

\begin{align}
c_{A\gamma}(x,y,z,t) = \sum_{\ell = 0}^{\ell=\infty} \sum_{m = 0}^{n=\infty}\sum_{n = 0}^{m=\infty}  & \exp\left( - ( \ell^2  +\tfrac{m^2}{d^2_{yy}}  + \tfrac{n^2}{d^2_{zz}} ) \frac{\pi^2}{L^2}D_{xx} t \right) \nonumber \\
 &\times B_{\ell} H_{m} S_{n} \cos(\ell \pi \tfrac{x}{L}) \cos(m \pi \tfrac{y}{L}) \cos(n \pi \tfrac{z}{L})
\end{align}

The constants $B_\ell, H_m$, and $S_n$ are found in the conventional manner using the initial condition and the orthogonally of the trigonometric functions.  Specifically, for $t=0$ we have
\begin{align}
\Phi(x,y,z) = \sum_{\ell = 0}^{\ell=\infty} \sum_{m = 0}^{n=\infty}\sum_{n = 0}^{m=\infty}   & B_{\ell} H_{m} S_{n} \cos(\ell \pi \tfrac{x}{L}) \cos(m \pi \tfrac{y}{L}) \cos(n \pi \tfrac{z}{L})
\end{align}
multiplying both sides of this expression by $ \cos(\ell' \pi \tfrac{x}{L}) \cos(m' \pi \tfrac{y}{L}) \cos(n' \pi \tfrac{z}{L})$ and integrating yields

\begin{align}
 \int_{x=0}^{x=L}  \int_{y=0}^{y=L}  \int_{z=0}^{z=L}  \Phi(x,y,z)& \cos(\ell' \pi \tfrac{x}{L}) \cos(m' \pi \tfrac{y}{L}) \cos(n' \pi \tfrac{z}{L})  dz\; dy\; dx  =  \nonumber  \\
&\sum_{\ell = 0}^{\ell=\infty} \sum_{m = 0}^{n=\infty}\sum_{n = 0}^{m=\infty}    
 \int_{x=0}^{x=L}  \int_{y=0}^{y=L}  \int_{z=0}^{z=L}
B_{\ell} H_{m} S_{n} \cos(\ell \pi \tfrac{x}{L})\cos(\ell' \pi \tfrac{x}{L}) \nonumber\\
& \times \cos(m \pi \tfrac{y}{L})  \cos(m' \pi \tfrac{y}{L})  \cos(n \pi \tfrac{z}{L})\cos(n' \pi \tfrac{z}{L})  dz\; dy\; dx
\end{align}
The right-hand side of this expression is non-zero only for the condition $\ell=\ell', m=m', n=n'$, for which we have

\begin{align}
 \int_{x=0}^{x=L}  \int_{y=0}^{y=L}  \int_{z=0}^{z=L}  \Phi(x,y,z)& \cos(\ell \pi \tfrac{x}{L}) \cos(m \pi \tfrac{y}{L}) \cos(n \pi \tfrac{z}{L})  dz\; dy\; dx  =  \nonumber  \\
& B_{\ell} H_{m} S_{n} \int_{x=0}^{x=L}  \int_{y=0}^{y=L}  \int_{z=0}^{z=L}
 \cos^2(\ell \pi \tfrac{x}{L})
  \cos^2(m \pi \tfrac{y}{L})   \cos^2(n \pi \tfrac{z}{L})  dz\; dy\; dx
\end{align}
The trigonometric integral on the right-hand side is easily verified to be $L^3/8$ when $\ell,m,n > 0$.  In general, the value of the integral is $L^3/(2^{N_{0}})$, where $N_{0}$ is the number of non-zero indexes (i.e., $N_{0}=0$, $1$, $2$ or $3$).  Thus, the value of $B_{\ell} H_{m} S_{n}$ is found from

\begin{equation}
 B_{\ell} H_{m} S_{n} = \frac{2^{N_0}}{L^3}  \int_{x=0}^{x=L}  \int_{y=0}^{y=L}  \int_{z=0}^{z=L}  \Phi(x,y,z) \cos(\ell \pi \tfrac{x}{L}) \cos(m \pi \tfrac{y}{L}) \cos(n \pi \tfrac{z}{L})  dz\; dy\; dx 
\end{equation}
Note here that the Fourier coefficients are not necessarily independently determinable.  For non-separable initial conditions, only the product of the coefficients can be determined; frequently a single symbol is used for this coefficient (e.g., $\bar{B}_{\ell m n} = B_{\ell} H_{m} S_{n}$).  The solution is given by

\begin{align}
c_{A\gamma}(x,y,z,t) =
 \sum_{\ell = 0}^{\ell=\infty} \sum_{m = 0}^{m=\infty} \sum_{n = 0}^{n=\infty} & \exp \left( -  [\ell^2  \tfrac{\pi^2}{L^2} D_{xx} + \tfrac{m^2}{d^2_{yy}}   \tfrac{\pi^2}{L^2} D_{xx}   +  \tfrac{n^2}{d^2_{zz}}  \tfrac{\pi^2}{L^2} D_{xx}] t \right) \nonumber \\
 &  \bar{B}_{\ell m n} \cos(\ell \pi \tfrac{x}{L})    \cos(m \pi \tfrac{y}{L})   \cos(n \pi \tfrac{z}{L})      
 \label{generalSol}
\end{align}
where the Fourier coefficients are evaluated from

\begin{equation}
 \bar{B}_{\ell m n} = \frac{2^{N_{0}}}{L^3}  \int_{x=0}^{x=L} \int_{y=0}^{y=L} \int_{z=0}^{z=L}  \Phi(x,y,z) \cos(\ell \pi \tfrac{x}{L})  \cos(m \pi \tfrac{y}{L}) \cos(n \pi \tfrac{z}{L}) \;dx \;dy  \;dz
 \label{fullFourier}
\end{equation}
For the non-separable initial conditions, it is necessary that the initial condition function be well-behaved enough such that a representation by a triple Fourier series is possible.  For the practical purposes discussed here, it is sufficient that the initial condition has meaning in at least a distributional sense.

There is a useful simplification for functions that {\it are multiplicative}  in the form $\Phi(x,y,z)=\Phi_x(x)\Phi_y(y)\Phi_z(z)$ (i.e., {\it separable}).  For this case, we have
\begin{equation}
 B_{\ell} H_{m} S_{n} = \frac{2^{N_{0}}}{L^3}  \int_{x=0}^{x=L}  \Phi_x(x) \cos(\ell \pi \tfrac{x}{L})\;dx 
 \int_{y=0}^{y=L}  \Phi_y(y) \cos(m \pi \tfrac{y}{L}) \;dy
 \int_{z=0}^{z=L}  \Phi_z(z)  \cos(n \pi \tfrac{z}{L})  \;dz
\end{equation}
and this allows us to determine $B_{\ell} , H_{m}$, and $S_{n}$ independently by

\begin{align}
B_{\ell} &= \frac{2^{ N_{0\ell} }}{L}  \int_{x=0}^{x=L}  \Phi_x(x) \cos(\ell \pi \tfrac{x}{L})\;dx \\
H_{m}&= \frac{2^{ N_{0 m} }}{L}   \int_{y=0}^{y=L}  \Phi_y(y) \cos(m \pi \tfrac{y}{L}) \;dy  \\
S_{n} & = \frac{2^{ N_{0 n} }}{L} \int_{z=0}^{z=L}  \Phi_z(z)  \cos(n \pi \tfrac{z}{L})  \;dz
\end{align}
where here, $N_{0\ell}, N_{0 m}$, and $N_{0 n}$ are equal to zero if the associated index ($\ell, m,$ or $n$) is zero, and equal to $1$ otherwise. For the problems considered here, we have imposed the condition that the total mass in the system is unity; this means that we can compute the zeroth coefficient immediately for the multiplicative initial condition directly by

\begin{align}
 B_0 &=  \frac{1}{L} \int_{x=0}^{x=L}  \int_{y=0}^{y=L}  \int_{z=0}^{z=L}   \Phi_x(x) dx =\frac{1}{L} \\
 H_0 &=  \frac{1}{L} \int_{x=0}^{x=L}  \int_{y=0}^{y=L}  \int_{z=0}^{z=L}   \Phi_y(y) dy =\frac{1}{L} \\
 S_0 &=  \frac{1}{L} \int_{x=0}^{x=L}  \int_{y=0}^{y=L}  \int_{z=0}^{z=L}   \Phi_z(z) dz =\frac{1}{L} 
\end{align}

Because the series involved is a cosine series, it is convenient to make the first term explicit, and represent the remaining terms as a sum from $1$ to $\infty$.  For the multiplicative initial condition, this yields a particularly simple expansion
\begin{align}
c_{A\gamma}(x,y,z,t) =&
\left[B_{0}+ \sum_{\ell = 1}^{\ell=\infty}    B_{\ell} \cos(\ell \pi \tfrac{x}{L})  \exp\left( -  \ell^2  \tfrac{\pi^2}{L^2} D_{xx} t \right)  \right]    \nonumber \\
& \times \left[  H_{0}+ \sum_{m = 1}^{m=\infty}  H_{m} \cos(m \pi \tfrac{y}{L})  \exp\left( - \tfrac{m^2}{d^2_{yy}}   \tfrac{\pi^2}{L^2} D_{xx} t \right)  \right]   \nonumber \\
&\times  \left[ S_{0}+ \sum_{n = 1}^{n=\infty}   S_{n}  \cos(n \pi \tfrac{z}{L})   \exp\left( -  \tfrac{n^2}{d^2_{zz}}  \tfrac{\pi^2}{L^2} D_{xx} t \right)  \right]  
\end{align}
%

\subsection{Moments}
%
In applications, it is often desirable to have explicit formulae for the spatial moments of the solution.  One reason for this is to provide a global (spatially) measure that can be compared with either experimental data or with numerical results to provide some sense of how much agreement there is with the analytical solution.  The zeroth, first, and second moments are usually three that are considered in applications.  The zeroeth moment (total mass) in each case is normalized to unity.  The first moment measures the center of mass of the solution as a function of time, and in the steady state it should be identically equal to ${\bf x}=(L/2, L/2, L/2)$ regardless of the initial condition.  The second centered moment is particularly useful for diffusion problems, because it relates to the spread of the solution around its first moment.  \\

{\noindent}{\it Zeroth moment}
\begin{equation}
m_0(t)  = \int_{x=0}^{x=L}  \int_{y=0}^{y=L}  \int_{z=0}^{z=L}  c_{A\gamma}(x,y,z,t) dz dy dx 
\end{equation}
For convenience, the zeroeth moment (total mass) for each case is normalized to unity, so that $m_0=1$.  To obtain a different value for the total mass, $m'_0$, one need only multiply the unit mass solution by the appropriate constant so that $c'(x,y,z,t) = m'_0 c(x,y,z,t)$, where $c'(x,y,z,t)$ is the concentration associated with the mass $m'_0$.  \\

There are three first moments for the concentration field, one for each of the three Cartesian axes.  The first moment is found by evaluating the following integrals. 

{\noindent}{\it First moment}
\begin{align}
m_{x}(t) &=  \int_{x=0}^{x=L}  \int_{y=0}^{y=L}  \int_{z=0}^{z=L}    x c_{A\gamma}(x,y,z,t) dz dy dx  \label{firstx} \\
m_{y}(t) &=  \int_{x=0}^{x=L}  \int_{y=0}^{y=L}  \int_{z=0}^{z=L}    y c_{A\gamma}(x,y,z,t) dz dy dx   \label{firsty} \\
m_{z}(t) &=  \int_{x=0}^{x=L}  \int_{y=0}^{y=L}  \int_{z=0}^{z=L}    z c_{A\gamma}(x,y,z,t) dz dy dx   \label{firstz}
\end{align}
For the case of multiplicative initial conditions, the first moments can be computed explicitly by
\begin{align}
m_{x}(t) &=  \frac{L^2}{2} \left[\frac{1}{L}+
 \frac{2}{\pi^2}\sum_{\ell=1}^{\ell=\infty} B_{\ell}   \frac{\left(-1  +(-1)^\ell \right)}{ \ell^2}    \exp\left( -  \ell^2  \tfrac{\pi^2}{L^2} D_{xx} t \right)    \right] \\
m_{y}(t) &=  \frac{L^2}{2} \left[\frac{1}{L}+
 \frac{2}{\pi^2}\sum_{m=1}^{m=\infty} H_{m}   \frac{\left(-1  +(-1)^m \right)}{ m^2}    \exp\left( - \tfrac{m^2}{d^2_{yy}}  \tfrac{\pi^2}{L^2} D_{xx} t \right)
    \right] \\
m_{z}(t) &=  \frac{L^2}{2} \left[\frac{1}{L}+
 \frac{2}{\pi^2}\sum_{n=1}^{n=\infty} S_{n}   \frac{\left(-1  +(-1)^n \right)}{ n^2}    \exp\left( -  \tfrac{n^2}{d^2_{zz}} \tfrac{\pi^2}{L^2} D_{xx} t \right)   \right]
\end{align}
For initial conditions that are symmetric about the point ${\bf x}=(L/2, L/2, L/2)$ for each of the three coordinate axes, the first moments are by definition $m_x = m_y=m_z = L/2$; thus they do not change in time.\\

The centered second moments are actually represented by a $3 \times 3$ tensor, of the form

\begin{equation}
 {\textbf{\sffamily \bfseries M}}=\left[ {\begin{array}{*{20}{c}}
  {{M_{xx}}}&{{M_{xy}}}&{{M_{xz}}} \\ 
  {{M_{yx}}}&{{M_{yy}}}&{{M_{yz}}} \\ 
  {{M_{zx}}}&{{M_{zy}}}&{{M_{zz}}} 
\end{array}} \right]
\end{equation}
In practice, however, it is generally the diagonal elements of the centered moment tensor are of primary interest.  The diagonal elements of the centered second moment tensor are given by the following.

{\noindent}{\it Centered second moment}
\begin{align}
M_{xx}(t) &=  \int_{x=0}^{x=L}  \int_{y=0}^{y=L}  \int_{z=0}^{z=L}    (x-\tfrac{L}{2})^2 c_{A\gamma}(x,y,z,t) dz dy dx   \\
M_{yy}(t) &=  \int_{x=0}^{x=L}  \int_{y=0}^{y=L}  \int_{z=0}^{z=L}    (y-\tfrac{L}{2})^2 c_{A\gamma}(x,y,z,t) dz dy dx   \\
M_{zz}(t) &=  \int_{x=0}^{x=L}  \int_{y=0}^{y=L}  \int_{z=0}^{z=L}    (z-\tfrac{L}{2})^2 c_{A\gamma}(x,y,z,t) dz dy dx  
\end{align}
For reference, the off-diagonal elements, $M_{\xi \eta}$, of the centered second moment tensor are defined analogously, e.g., 

\begin{equation}
M_{xy} =  \int_{x=0}^{x=L}  \int_{y=0}^{y=L}  \int_{z=0}^{z=L}    (x-\tfrac{L}{2})(y-\tfrac{L}{2}) c_{A\gamma}(x,y,z,t) dz dy dx
\end{equation} 
however, the off-diagonal elements will not be considered further here.

For the case of multiplicative initial conditions, an explicit representation of the second moments can be determined by evaluating the integrals above.  This leads to the closed-form expressions
\begin{align}
M_{xx}(t) &=  \frac{L^3}{12} \left[ \frac{1}{L}+
 \frac{12}{\pi^2}\sum_{\ell=1}^{\ell=\infty} B_{\ell}   \frac{\left(1  +(-1)^\ell \right)}{ \ell^2}    \exp\left( -  \ell^2  \tfrac{\pi^2}{L^2} D_{xx} t \right)   \right] \\
M_{yy}(t) &=  \frac{L^3}{12} \left[ \frac{1}{L}+
 \frac{12}{\pi^2}\sum_{m=1}^{m=\infty} H_{m}   \frac{\left(1  +(-1)^m \right)}{ m^2}     \exp\left( - \tfrac{m^2}{d^2_{yy}}  \tfrac{\pi^2}{L^2} D_{xx} t \right)    \right] \\
M_{zz}(t) &=  \frac{L^3}{12} \left[\frac{1}{L}+
 \frac{12}{\pi^2}\sum_{n=1}^{n=\infty} S_{n}   \frac{\left(1  +(-1)^n \right)}{ n^2}    \exp\left( -  \tfrac{n^2}{d^2_{zz}} \tfrac{\pi^2}{L^2} D_{xx} t \right)     \right]
\end{align}
It is clear from these expressions that at arbitrarily large times, these moments all tend toward the value $L^2/12$, which is the correct value for the centered moment of a uniform cube.

\section{Explicit Series Solutions for Particular Initial Conditions}

In this section, the series solution for the four initial conditions described previously are computed.  The first (when it is different from the value $\tfrac{L}{2}$) and second moments are also computed as series.  To provide a more intuitive feel for the solutions, isosurface plots of the concentration field are presented.  For presentation of the plots, we define the diffusive time scales in each of the cardinal directions by

\begin{align}
T_x^* =& \frac{(\tfrac{L}{2})^2}{{D_{xx}}} \label{tstarx}  \\
T_y^* =& \frac{(\tfrac{L}{2})^2}{d^2_{yy}{D_{xx}}} \label{tstary}\\
T_z^* =& \frac{(\tfrac{L}{2})^2}{d^2_{zz}{D_{xx}}}\label{tstarz}
\end{align}
These represent, very approximately, the time it takes for diffusive gradients in the each direction to move a substantial amount of mass to the boundary.  Intuitively, gradients in the each direction should approach zero after just a few multiples of the diffusive time scale for that direction.  For the plots generated in this work, we will use the set of parameters defined in Table \ref{table-parameters}; these parameters represent physically reasonable values for diffusion of dilute aqueous species in fluids or gels.  For other physical systems (e.g., heat transfer, diffusion in solids), these parameters would obviously have different values.  However, there is substantial utility in providing examples that correspond physically to some system, as we have done here.   As seen in Table \ref{table-parameters}, $T^*_x$ is the smallest time scale.  Thus, for plotting we have adopted the following normalized time scale

\begin{equation}
t^* = \frac{t}{T^*_x}
\end{equation}

To provide additional confidence in the series solutions have been computed, each using only the first $20$ terms of the series, and compared with the results of a previously verified numerical code.  The code and convergence is described in detail in \S \ref{computations}.  

\begin{table}
  \caption{Parameters used in the example computations.}
  \begin{center}
  \begin{tabular}{c c c c}
 
  Parameter 	& Description    &	Value  & Units\\[5pt]
  \hline \\
  $L$ 		& Domain Size   &	$0.01$ & [m] \\
  $m_0$             &Total mass of diffusing solute   $1.0$   [$\mu g$]  \\
  $c_\infty=\tfrac{m_0}{L^3}$         &Steady-state solute concentration   & $1\times10^{6}$  &[$\mu g/m^3$] \\
  $M_\infty=\tfrac{L^2}{12}$           &Steady-state second moment  &  $8.33\times10^{-6}$  &[$m^2$]  \\
  ${D_{xx}}$      &Diffusion coefficient, $x-$direction  &	$1\times10^{-9}$ & [$m^2$/s] \\
  $d^2_{yy}$ 	& Anisotorpy ratio, $y-$direction  &	${2}$ & [-]\\
  $d^2_{zz}$ 	& Anisotorpy ratio, $z-$direction	&4 & [-]\\
  $T^*_x$            & Diffusion time scale, $x-$direction   &  25000  & [s] \\
  $T^*_y$            &  Diffusion time scale, $y-$direction &   50000  & [s] \\
  $T^*_z$            &   Diffusion time scale, $z-$direction &  100000  & [s] \\
  $a$ 		& Size parameter for initial condition, Case 2 &	0.5$L$ & [m]\\
  $\sigma_x$ 	& Variance for initial condition, Case 3&	0.1$L$ & [m]\\
  $\kappa_y$    & Parameter for planar initial condition, Case 4 &  20 &  [-] \\
    $\kappa_z$    & Parameter for planar initial condition, Case 4 &  40 &  [-] \\
  \hline \\
  \end{tabular}
   \label{table-parameters} 
  \end{center}
\end{table}

\subsection{Case 1}

For Case 1, the initial condition is given by the delta function.  The delta function provides a good case for assessing how well codes conserve the dependent variable for challenging initial conditions, and are also useful for verifying particle tracking codes \citep{xu2010verification}.  The initial condition for this case is specified by

\begin{align}
&I.C.& c({\bf x},0)&= \delta({\bf x}-{\bf L}/2) = \delta(x-L/2)\delta(y-L/2)\delta(z-L/2)
\end{align}

\begin{figure}
  \centering
 \includegraphics[scale=0.15]{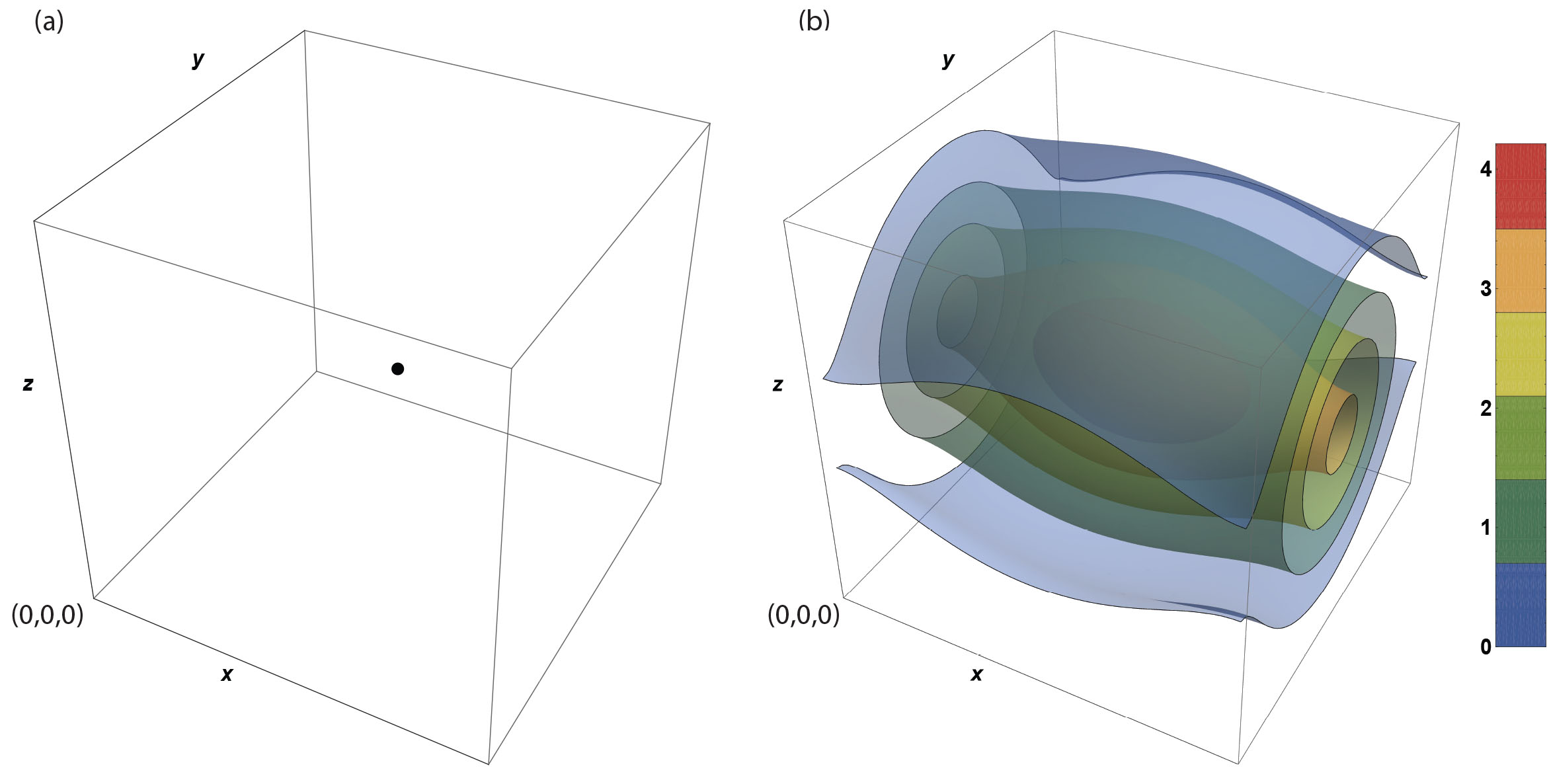}
  \caption{Concentration isosurfaces for the $\delta$-function initial condition (Case 1);  concentration is normalized, $\tfrac{c}{c_\infty}$.  (a) Although formally not plottable, the initial condition is represented as a point distribution with finite mass. (b) The concentration field at $t^*=\tfrac{1}{4} T^*_x$ using the first 20 terms of the series solution.  }
    \label{delta_concentration}
\end{figure}

The values of the Fourier series constants are found by the conventional arguments using the orthogonality principles of trigonometric series.  This evaluation leads to

\begin{align}
B_{0} &=H_{0}=S_{0}=  \frac{1}{L} \\
B_{\ell} &= \frac{2}{L} \int_0^L   \delta({ x} -{L}/2) cos(\ell \pi \tfrac{x}{L})  dx  = \frac{2}{L}cos( \tfrac{\ell \pi}{2}) , ~for ~ \ell  > 0 \\
H_{m} &= \frac{2}{L} \int_0^L   \delta({ y} -{L}/2) cos(m \pi \tfrac{y}{L})  dy  = \frac{2}{L}cos( \tfrac{m \pi}{2}) , ~for ~ m  > 0 \\
S_{n} &= \frac{2}{L} \int_0^L  \delta({ z} -{L}/2) cos(n \pi \tfrac{z}{L})  dz  = \frac{2}{L}cos( \tfrac{n \pi}{2}) , ~for ~ n  > 0
\end{align}
Hence, the solution for the concentration field  can be given explicitly by

\begin{align}
c_{A\gamma}(x,y,z,t) =&
\frac{1}{L^3}\left[1+ 2\sum_{\ell = 1}^{\ell=\infty}    \cos(\tfrac{\ell \pi}{2}) \cos(\ell \pi \tfrac{x}{L})  \exp\left( -  \ell^2  \tfrac{\pi^2}{L^2} D_{xx} t \right)  \right]    \nonumber \\
& \times \left[1+ 2\sum_{m = 1}^{m=\infty}   \cos(\tfrac{m \pi}{2}) \cos(m \pi \tfrac{y}{L})  \exp\left( - \tfrac{m^2}{d^2_{yy}}  \tfrac{\pi^2}{L^2} D_{xx} t \right)  \right]   \nonumber \\
&\times  \left[1+ 2\sum_{n = 1}^{n=\infty}    \cos(\tfrac{n \pi}{2})  \cos(n \pi \tfrac{z}{L})   \exp\left( -  \tfrac{n^2}{d^2_{zz}} \tfrac{\pi^2}{L^2} D_{xx} t \right)  \right]  
\end{align}
Isosurface plots of the normalized concentration appear in Fig.~\ref{delta_concentration}.\footnote{Mathematica scripts for computing each of the four concentration field solutions and associated moments as series will be included as ancillary materials in the final submission of this manuscript.}

The moments for this solution are found as follows.\\

{\noindent}{\it Zeroth moment}
\begin{equation}
m_0(t)  = \int_{x=0}^{x=L}  \int_{y=0}^{y=L}  \int_{z=0}^{z=L}  c_{A\gamma}(x,y,z,t) dz dy dx = 1
\end{equation}

{\noindent}{\it First moment}
\begin{equation}
m_x(t)  = m_y(t) = m_z(t) = \frac{L}{2}
\end{equation}
The zeroeth and first moments are constant for all cases considered, except Case 4 (where the first moment is transient).  Where the moments are constant, they are not plotted.

{\noindent}{\it Centered second moments}
\begin{align}
M_{xx}(t) &=  \int_{x=0}^{x=L}  \int_{y=0}^{y=L}  \int_{z=0}^{z=L}    (x-\tfrac{L}{2})^2 c_{A\gamma}(x,y,z,t) dz dy dx  \nonumber \\
&=\frac{L^2}{12}  \left[1+
\sum_{\ell=1}^{\ell=\infty}\frac{24}{{\ell}^2 \pi ^2}\exp\left(-\ell^2 \tfrac{\pi ^2 }{L^2}D_{xx} t\right) \cos \left(\tfrac{{\ell} \pi }{2}\right) \left(1  +(-1)^\ell \right)\right] \\
M_{yy}(t) &=  \int_{x=0}^{x=L}  \int_{y=0}^{y=L}  \int_{z=0}^{z=L}    (y-\tfrac{L}{2})^2 c_{A\gamma}(x,y,z,t) dz dy dx  \nonumber \\
&=\frac{L^2}{12}  \left[ 1+
\sum_{m=1}^{m=\infty}\frac{24}{{m}^2 \pi ^2}\exp\left(-\tfrac{m^2}{d^2_{yy}}\tfrac{\pi ^2 }{L^2}D_{xx} t\right) \cos \left(\tfrac{{m} \pi }{2}\right) \left(1  +(-1)^m \right)\right] \\
M_{zz}(t) &=  \int_{x=0}^{x=L}  \int_{y=0}^{y=L}  \int_{z=0}^{z=L}    (z-\tfrac{L}{2})^2 c_{A\gamma}(x,y,z,t) dz dy dx  \nonumber \\
&=\frac{L^2}{12}  \left[1+
\sum_{n=1}^{n=\infty}\frac{24}{{n}^2 \pi ^2}\exp\left(-\tfrac{n^2}{d^2_{zz}}\tfrac{\pi ^2 }{L^2}D_{xx} t\right) \cos \left(\tfrac{{n} \pi }{2}\right) \left(1  +(-1)^n \right)\right] 
\end{align}

\begin{figure}
  \centering
 \includegraphics[scale=0.8]{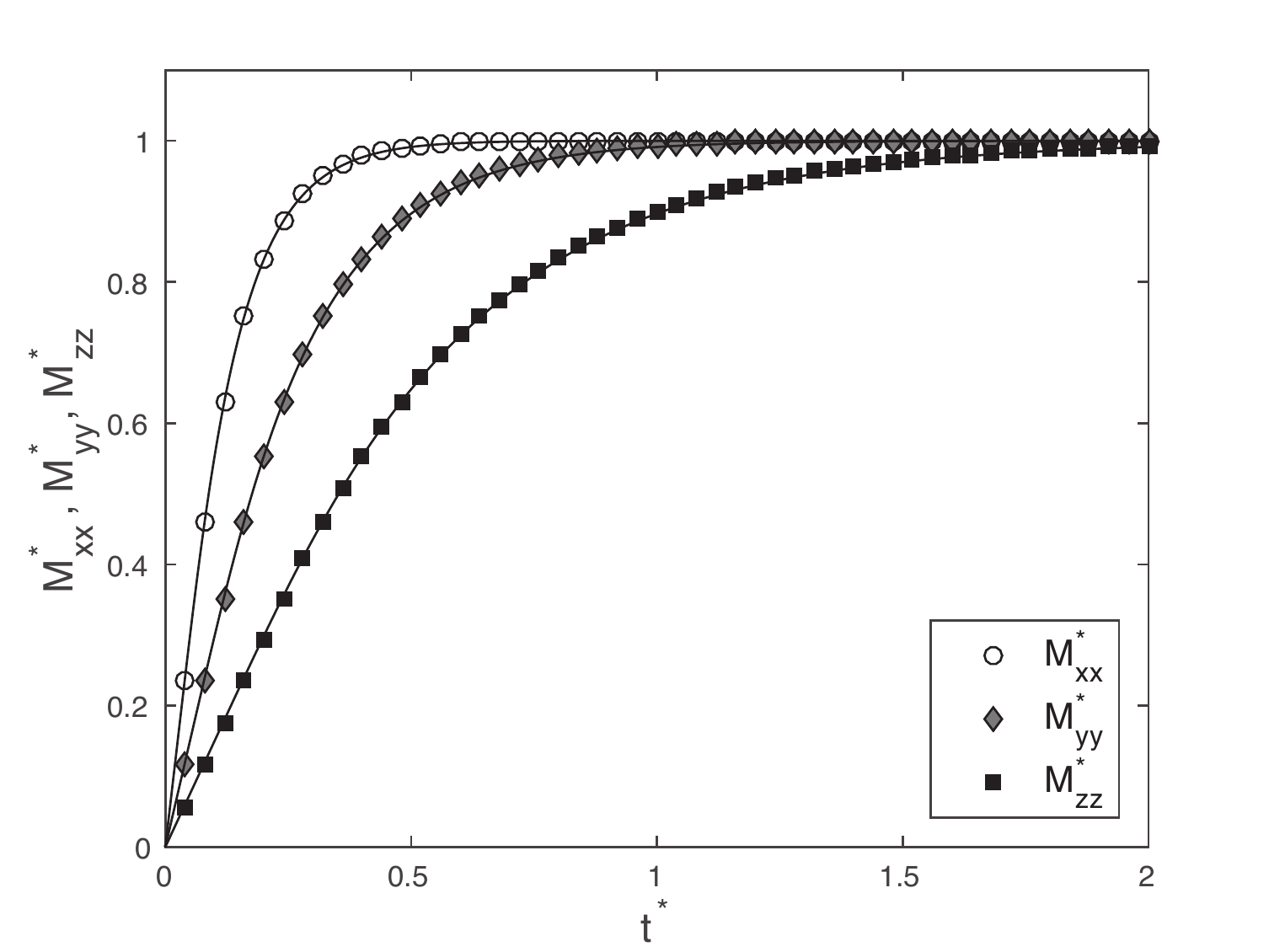}
  \caption{Centered second moments for the $\delta$-function initial condition (Case 1) via series solution (solid lines), and numerical computations  (symbols).}
    \label{delta}
\end{figure}

The centered second moments for Case 1 are plotted in Fig.~\ref{delta}.  For these plots, we have normalized the moments by their steady-state value, i.e., 

\begin{equation}
{M_\infty } = \mathop {\lim }\limits_{t \to \infty } [{M_{xx}}(t)] = \mathop {\lim }\limits_{t \to \infty } [{M_{yy}}(t)] = \mathop {\lim }\limits_{t \to \infty } [{M_{zz}}(t)] = \frac{{{L^2}}}{{12}}
\end{equation}
For plotting purposes, we have defined the following normalized variables

\begin{align}
M^*_{xx}=& \frac{M_{xx}}{M_\infty }   \\
M^*_{yy}=& \frac{M_{yy}}{M_\infty }  \\
M^*_{zz}=& \frac{M_{zz}}{M_\infty }  
\end{align}

\subsection{Case 2}

%
For case 2, the initial condition is given by a uniform concentration in a cube of side $a$ centered in the domain.  This initial condition poses the typical challenges that discontinuities create in numerical approximation.  The solution is also complex enough that it satisfies the guidelines specified by \citet{roache2002code} for defining {\it good} analytical solutions for code verification.  To define the initial condition, we first define the following functions in terms of unit Heaviside step functions.

\begin{align}
\Phi_x(x) =& H[x-(\tfrac{L}{2}-\tfrac{a}{2})]-H[x-(\tfrac{L}{2}+\tfrac{a}{2})]   \\
\Phi_y(y) =& H[y-(\tfrac{L}{2}-\tfrac{a}{2})]-H[y-(\tfrac{L}{2}+\tfrac{a}{2})]   \\
\Phi_z(z) =& H[z-(\tfrac{L}{2}-\tfrac{a}{2})]-H[z-(\tfrac{L}{2}+\tfrac{a}{2})]   
\end{align}
Where $H(x)$ is the unit Heaviside step function.  The function, $\Phi({\bf x})$, is given by
Defining the function $\Phi({\bf x})$ by

\begin{equation}
\Phi({\bf x}) = \Phi_x(x) \Phi_y(y) \Phi_z(z)
\end{equation}
The initial condition for a 3-dimensional step function centered in the cubic domain is now specified by

\begin{align}
&I.C.& c({\bf x},0)&=  \Phi({\bf x})
\end{align}

The values of the Fourier series constants are found in the conventional manner. 

\begin{align}
B_{0} &=H_{0}=S_{0}=  \frac{1}{L} \\
B_{\ell} &= \frac{2}{L} \int_0^L  \Phi(x,y,z) \cos(\ell \pi \tfrac{x}{L})  dx  = \frac{4}{L} \int_{\frac{L}{2}-\frac{a}{2}}^{\frac{L}{2}+\frac{a}{2}} \cos(\ell \pi \tfrac{x}{L})   dx = \nonumber\\
&\frac{4 \cos \left(\frac{\pi  \ell}{2}\right) \sin \left(\frac{\pi  a \ell}{2 L}\right)}{\pi a \ell}, ~for ~ \ell  > 0\\
H_{m} &=  \frac{2}{L} \int_0^L  \Phi(x,y,z) \cos(m \pi \tfrac{y}{L})  dy  = \frac{2}{L} \int_{\frac{L}{2}-\frac{a}{2}}^{\frac{L}{2}+\frac{a}{2}} \cos(m \pi \tfrac{y}{L})   dy = \nonumber\\
&\frac{4 \cos \left(\frac{\pi  m}{2}\right) \sin \left(\frac{\pi  a m}{2 L}\right)}{\pi a m}, ~for ~ m  > 0 \\
S_{n} &=    \frac{2}{L} \int_0^L  \Phi(x,y,z) \cos(n \pi \tfrac{z}{L})  dz  = \frac{2}{L} \int_{\frac{L}{2}-\frac{a}{2}}^{\frac{L}{2}+\frac{a}{2}} \cos(n \pi \tfrac{z}{L})   dz = \nonumber\\
&\frac{4 \cos \left(\frac{\pi  n}{2}\right) \sin \left(\frac{\pi  a n}{2 L}\right)}{\pi a n}, ~for ~ n  > 0
\end{align}

\begin{figure}
  \centering
 \includegraphics[scale=0.15]{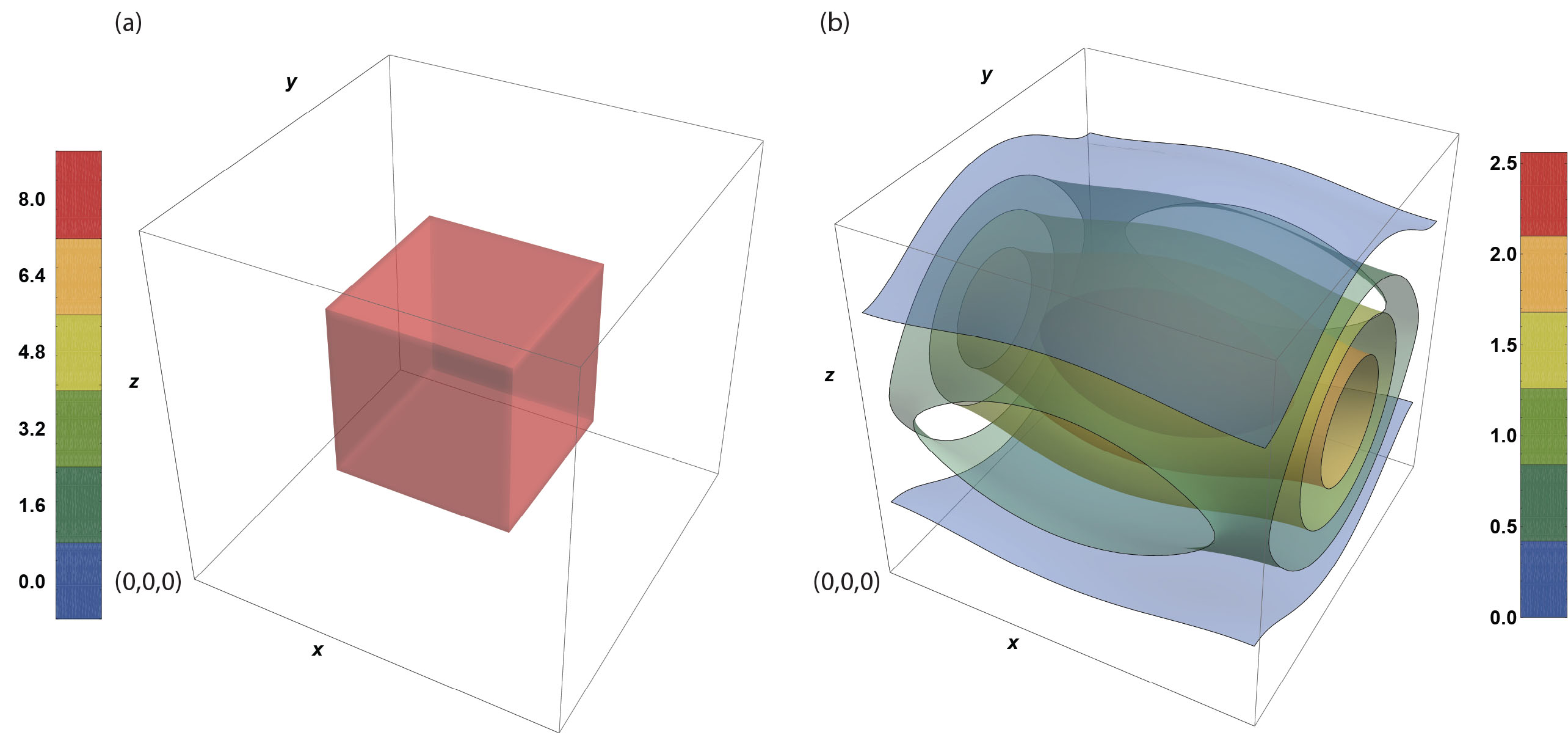}
  \caption{Concentration isosurfaces for the step function initial condition (Case 2);  concentration is normalized, $\tfrac{c}{c_\infty}$.  (a) The initial condition. (b) The concentration field at $t^*=\tfrac{1}{4} T^*_x$ using the first 20 terms of the series solution.  }
    \label{step_concentration}
\end{figure}

The explicit series expression for the concentration field is

\begin{align}
c_{A\gamma}(x,y,z,t) =&
\frac{1}{L^3}
\left[1+ \sum_{\ell = 1}^{\ell=\infty}   
\frac{4 L \cos \left(\frac{\pi  \ell}{2}\right) \sin \left(\frac{\pi  a \ell}{2 L}\right)}{\pi  a \ell} \cos(\ell \pi \tfrac{x}{L})
\exp\left( -  \ell^2  \tfrac{\pi^2}{L^2} D_{xx} t \right) 
   \right]    \nonumber \\
& \times \left[  1+ \sum_{m = 1}^{m=\infty}  
\frac{4 L \cos \left(\frac{\pi  m}{2}\right) \sin \left(\frac{\pi  a m}{2 L}\right)}{\pi  a m} \cos(m \pi \tfrac{y}{L})
 \exp\left( - \tfrac{m^2}{d^2_{yy}}   \tfrac{\pi^2}{L^2} D_{xx} t \right)
  \right]   \nonumber \\
&\times  \left[ 1+ \sum_{n = 1}^{n=\infty} 
\frac{4 L \cos \left(\frac{\pi  n}{2}\right) \sin \left(\frac{\pi  a n}{2 L}\right)}{\pi  a n} \cos(n \pi \tfrac{z}{L})
 \exp\left( - \tfrac{n^2}{d^2_{zz}}   \tfrac{\pi^2}{L^2} D_{xx} t \right)
    \right]  
\end{align}
Isosurface plots of the normalized concentration appear in Fig.~\ref{step_concentration}.

The moments for this initial condition are given by the following.\\

{\noindent}{\it Zeroth moment}
\begin{equation}
m_0(t)  = \int_{x=0}^{x=L}  \int_{y=0}^{y=L}  \int_{z=0}^{z=L}  c_{A\gamma}(x,y,z,t) dz dy dx = 1
\end{equation}

{\noindent}{\it First moment}
\begin{equation}
m_x(t)  = m_y(t) = m_z(t) = \frac{L}{2}
\end{equation}
As with Case 1, these moments are constant in time.  The centered second moments are, however, transient.  Explicit series solutions for these are given by\\

\begin{figure}
  \centering
 \includegraphics[scale=0.8]{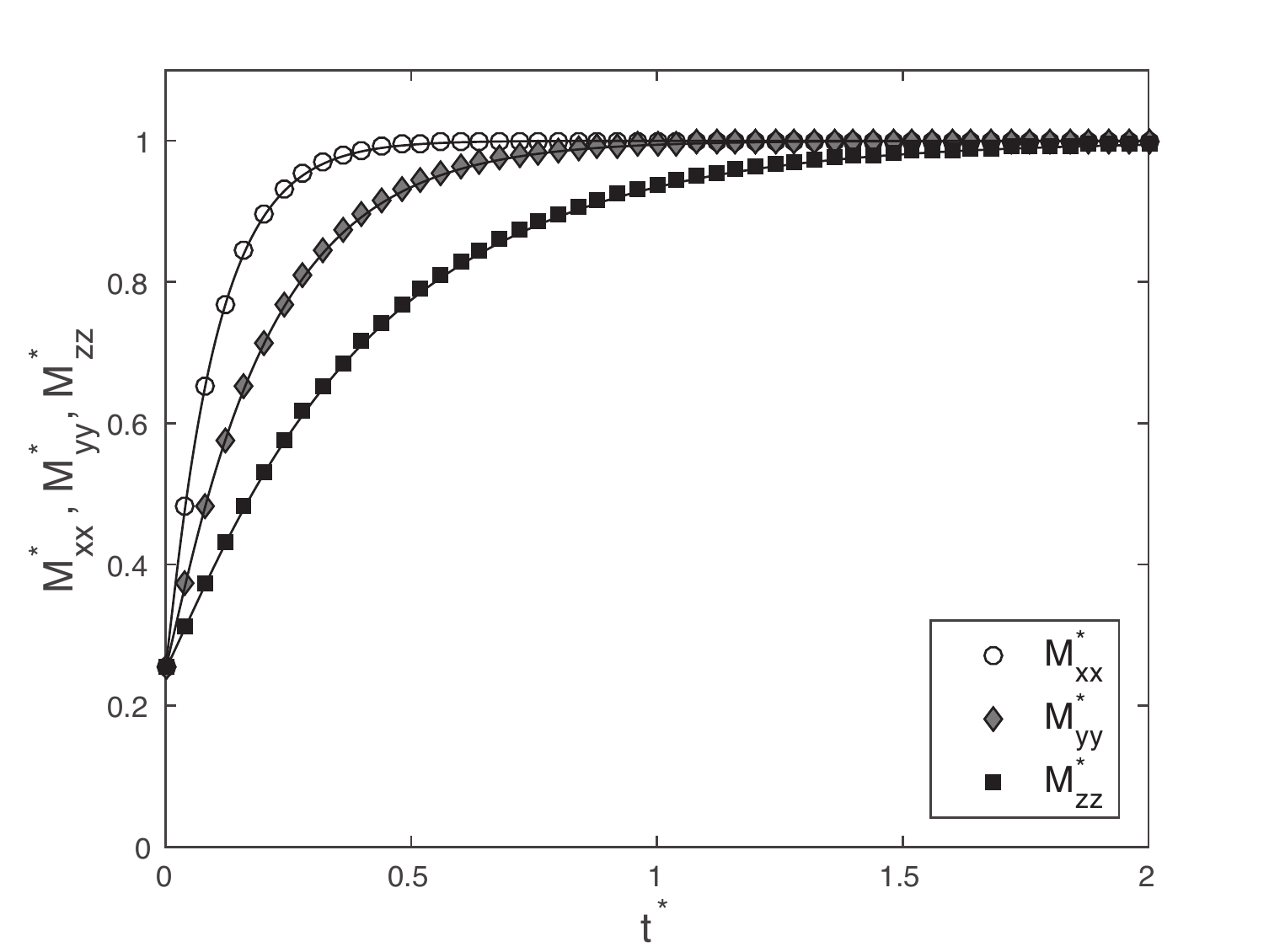}
  \caption{Centered second moments for the step initial condition (Case 2) via series solution (solid lines), and numerical computations  (symbols).}
    \label{step}
\end{figure}

{\noindent}{\it Centered second moments}
\begin{align}
M_{xx}(t) &=\frac{L^2}{12} \left[1+
 \frac{12}{\pi^2}\sum_{\ell=1}^{\ell=\infty} \frac{4 L \cos \left(\frac{\pi  \ell}{2}\right) \sin \left(\frac{\pi  a \ell}{2 L}\right)}{\pi a \ell} \frac{\left(1  +(-1)^\ell \right)}{ \ell^2}   \exp\left( -  \ell^2  \tfrac{\pi^2}{L^2} D_{xx} t \right)
   \right] \\
M_{yy}(t) &=\frac{L^2}{12} \left[1+
 \frac{12}{\pi^2}\sum_{m=1}^{m=\infty} \frac{4 L \cos \left(\frac{\pi  m}{2}\right) \sin \left(\frac{\pi  a m}{2 L}\right)}{\pi a m}   \frac{\left(1  +(-1)^m \right)}{ m^2}     \exp\left( - \tfrac{m^2}{d^2_{yy}}  \tfrac{\pi^2}{L^2} D_{xx} t \right)    \right] \\
M_{zz}(t) &=\frac{L^2}{12} \left[1+
 \frac{12}{\pi^2}\sum_{n=1}^{n=\infty} \frac{4 L \cos \left(\frac{\pi  n}{2}\right) \sin \left(\frac{\pi  a n}{2 L}\right)}{\pi  a n}   \frac{\left(1  +(-1)^n \right)}{ n^2}     \exp\left( -  \tfrac{n^2}{d^2_{zz}} \tfrac{\pi^2}{L^2} D_{xx} t \right)     \right]
\end{align}
In Fig.~\ref{step}, the normalized second moments are plotted as a function of time.

\subsection{Case 3}

\begin{figure}
  \centering
 \includegraphics[scale=0.15]{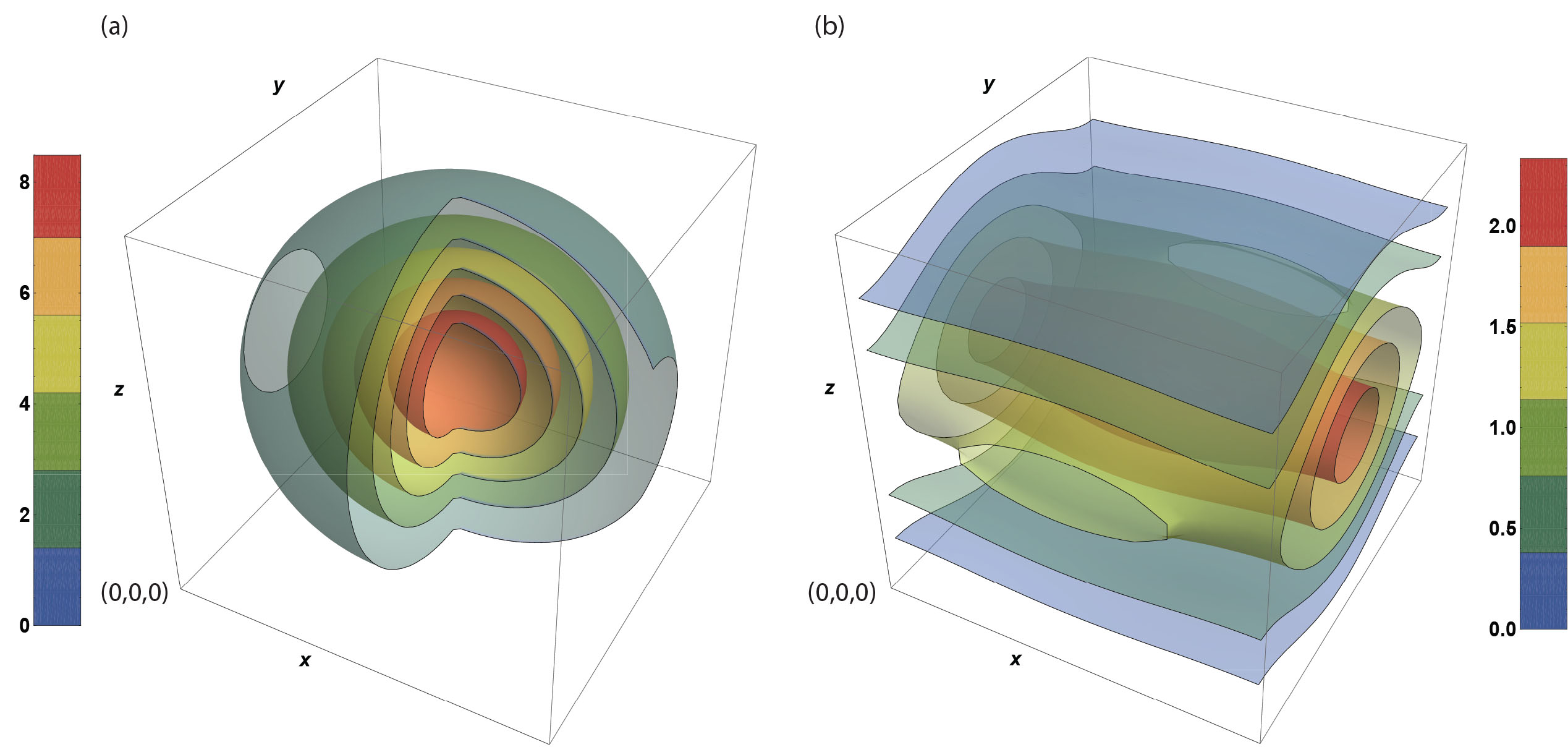}
  \caption{Concentration isosurfaces for the truncated Gaussian initial condition (Case 3);  concentration is normalized, $\tfrac{c}{c_\infty}$.  (a) The initial condition. (b) The concentration field at $t^*=\tfrac{1}{4} T^*_x$ using the first 20 terms of the series solution.  }
    \label{gaussian_concentration}
\end{figure}

For case 3, the initial condition is given by a truncated Gaussian function.  This function is a Gaussian function, centered at $x=y=z=L/2$ which has been truncated at the boundaries of the domain, and normalized so that the total integrated mass in the domain is unity.   Because the function is analytically smooth, numerical approximations to derivatives of it are well-behaved; however, it is complex enough to provide a challenging validation benchmark \citep{roache2002code}.  This initial condition function is specified explicitly by

\begin{align}
&I.C.& c({\bf x},0)&= 
\frac{ \exp \left(-\frac{\left(x-\frac{L}{2}\right)^2}{2
   {\sigma^2_x}}-\frac{
   \left(y-\frac{L}{2}\right)^2}{2 {\sigma^2_x}/{d^2_{yy}}}-\frac{
   \left(z-\frac{L}{2}\right)^2}{2 {\sigma^2_x}/{d^2_{zz}}}\right)}
{\sqrt{2  \pi \sigma^2_x }  \sqrt{2  \pi \sigma^2_x/d^2_{yy} }  \sqrt{2  \pi \sigma^2_x/d^2_{zz} }\,\,
 \text{erf}\left(\frac{L}{2 \sqrt{2} {{\sigma_x}}}\right) 
\text{erf}\left(\frac{ L}{2\sqrt{2} {(\sigma_x}/d_{yy})}\right)
\text{erf}\left(\frac{ L}{2\sqrt{2} {(\sigma_x}/d_{zz})}\right)
}
\end{align}

\begin{align}
B_0 &= H_0 = S_0 =\frac{1}{L} \\
B_{\ell} &= \frac{2}{L} \frac{  \Re\left[
   \text{erf}\left(\frac{L^2+2 i \pi  {\ell} \sigma_x^2}{2
   \sqrt{2} L {\sigma_x}}\right)\right]}
   { \text{erf}\left(\frac{L}{2 \sqrt{2}
   {\sigma_x}}\right)}
    \cos \left(\tfrac{\pi  {\ell}}{2}\right) 
  \exp\left(-{\ell}^2\tfrac{\pi ^2}{L^2}  \tfrac {\sigma_x^2}{2 } \right)  , ~for ~ \ell  > 0  \\
 H_m & = \frac{2}{L} \frac{  \Re\left[
   \text{erf}\left(\frac{{d^2_{yy}} L^2+2 i \pi  m {\sigma^2_x}}{2 \sqrt{2}
   L {{d_{yy}} {\sigma_x}}}\right)\right]}
   {\text{erf}\left(\frac{d_{yy} L}{2 \sqrt{2}
   {\sigma_x}}\right)}
    \cos \left(\tfrac{\pi  m}{2}\right) 
 \exp\left(-\tfrac{m^2}{d_{yy}^2}\tfrac{\pi ^2}{L^2}  \tfrac {\sigma_x^2}{2 } \right)
  , ~for ~ m  > 0   \\
S_n &=\frac{2}{L}  \frac{ 
 \Re\left[\text{erf}\left(\frac{{d^2_{zz}} L^2+2 i \pi  m \sigma_x^2}{2 \sqrt{2}
   L {{d_{zz}} \sigma_x}}\right)\right]}
   { \text{erf}\left(\frac{d_{zz} L}{2 \sqrt{2}
   {\sigma_x}}\right)}  
   \cos \left(\tfrac{\pi  n}{2}\right) 
 \exp\left(-\tfrac{n^2}{d_{zz}^2}\tfrac{\pi ^2}{L^2}  \tfrac {\sigma_x^2}{2 } \right)
 , ~for ~ n  > 0 
\end{align}
Note that here  at the limit as $\sigma_x \rightarrow 0$, this expression reduces to that for the delta initial condition.  Performing the integrations required to evaluate the Fourier coefficients, the closed-form solution for the concentration field is

\begin{align}
&c_{A\gamma}(x,y,z,t) =\\
&
\frac{1}{L^3}\left[1+ \sum_{\ell = 1}^{\ell=\infty}    2 \frac{  \Re\left[
   \text{erf}\left(\frac{L^2+2 i \pi  {\ell} \sigma_x^2}{2
   \sqrt{2} L {\sigma_x}}\right)\right]}
   { \text{erf}\left(\frac{L}{2 \sqrt{2}
   {\sigma_x}}\right)}
    \cos \left(\tfrac{\pi  {\ell}}{2}\right) 
    \exp\left(-{\ell}^2\tfrac{\pi ^2}{L^2}  \tfrac {\sigma_x^2}{2 } \right)
   \cos(\ell \pi \tfrac{x}{L})  \exp\left( -  \ell^2  \tfrac{\pi^2}{L^2} D_{xx} t \right)  \right]    \nonumber \\
& \times \left[  1+ \sum_{m = 1}^{m=\infty}  2 \frac{  \Re\left[
   \text{erf}\left(\frac{{d^2_{yy}} L^2+2 i \pi  m {\sigma^2_x}}{2 \sqrt{2}
   L {{d_{yy}} {\sigma_x}}}\right)\right]}
   {\text{erf}\left(\frac{d_{yy} L}{2 \sqrt{2}
   {\sigma_x}}\right)}
    \cos \left(\tfrac{\pi  m}{2}\right) 
 \exp\left(-\tfrac{m^2}{d_{yy}^2}\tfrac{\pi ^2}{L^2}  \tfrac {\sigma_x^2}{2 } \right)
   \cos(m \pi \tfrac{y}{L})  \exp\left( - \tfrac{m^2}{d^2_{yy}}   \tfrac{\pi^2}{L^2} D_{xx} t \right)  \right]   \nonumber \\
&\times  \left[ 1+ \sum_{n = 1}^{n=\infty}   2 \frac{ 
 \Re\left[\text{erf}\left(\frac{{d^2_{zz}} L^2+2 i \pi  n \sigma_x^2}{2 \sqrt{2}
   L {{d_{zz}} \sigma_x}}\right)\right]}
   { \text{erf}\left(\frac{d_{zz} L}{2 \sqrt{2}
   {\sigma_x}}\right)}  
   \cos \left(\tfrac{\pi  n}{2}\right) 
\exp\left(-\tfrac{n^2}{d_{zz}^2}\tfrac{\pi ^2}{L^2}  \tfrac {\sigma_x^2}{2 } \right)
   \cos(n \pi \tfrac{z}{L})   \exp\left( -  \tfrac{n^2}{d^2_{zz}}  \tfrac{\pi^2}{L^2} D_{xx} t \right)  \right]  
\end{align}\\
Isosurface plots of the normalized concentration appear in Fig.~\ref{gaussian_concentration}.

The moments for this initial condition are given by the following.\\

{\noindent}{\it Zeroth moment}
\begin{equation}
m_0(t)  = \int_{x=0}^{x=L}  \int_{y=0}^{y=L}  \int_{z=0}^{z=L}  c_{A\gamma}(x,y,z,t) dz dy dx = 1
\end{equation}

{\noindent}{\it First moment}
\begin{equation}
m_x(t)  = m_y(t) = m_z(t) = \frac{L}{2}
\end{equation}
As before, the first two moments are constant in time.  The centered second moments are transient, and specified by

\begin{figure}
  \centering
 \includegraphics[scale=0.8]{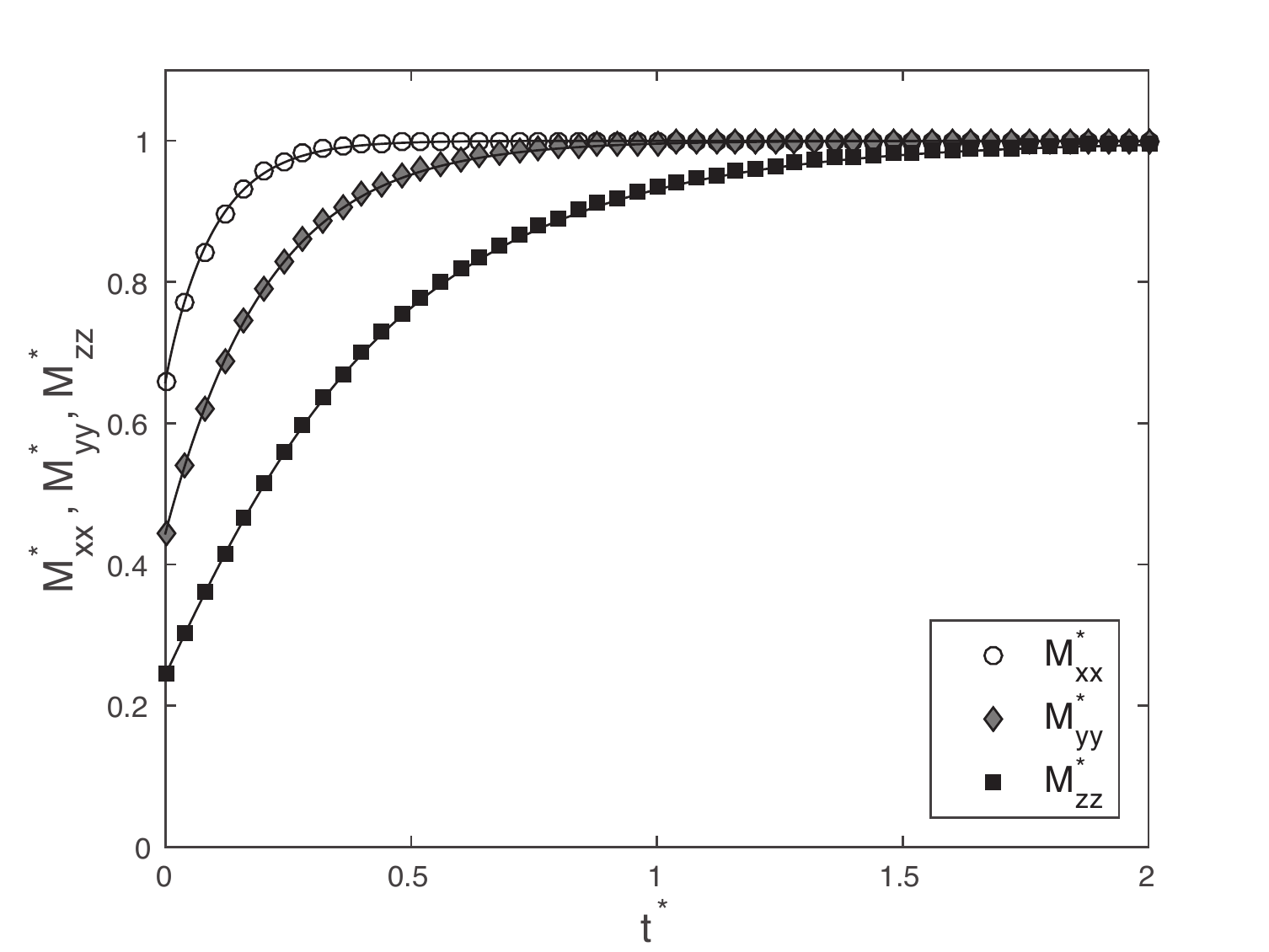}
  \caption{Centered second moments for the truncated Gaussian initial condition (Case 3) via series solution (solid lines), and numerical computations  (symbols).}
      \label{gaussian}
\end{figure}

{\noindent}{\it Centered second moments}
\begin{align}
M_{xx}(t) &= \frac{L^2}{12} \Bigg[1+
 \frac{24}{\pi^2}\sum_{\ell=1}^{\ell=\infty}  \frac{  \Re\left[
   \text{erf}\left(\frac{L^2+2 i \pi  {\ell} \sigma_x^2}{2
   \sqrt{2} L {\sigma_x}}\right)\right]}
   { \text{erf}\left(\frac{L}{2 \sqrt{2}
   {\sigma_x}}\right)}
    \cos \left(\tfrac{\pi  {\ell}}{2}\right) \nonumber \\
&\times \exp\left(-{\ell}^2\tfrac{\pi ^2}{L^2}  \tfrac {\sigma_x^2}{2 } \right)   
   \frac{\left(1  +(-1)^\ell \right)}{ \ell^2}     \exp\left( -  \ell^2  \tfrac{\pi^2}{L^2} D_{xx} t \right) \Bigg] \\
M_{yy}(t) &=\frac{L^2}{12} \Bigg[1+
 \frac{24}{\pi^2}\sum_{m=1}^{m=\infty} \frac{  \Re\left[
   \text{erf}\left(\frac{{d^2_{yy}} L^2+2 i \pi  m {\sigma^2_x}}{2 \sqrt{2}
   L {{d_{yy}} {\sigma_x}}}\right)\right]}
   {\text{erf}\left(\frac{d_{yy}L}{2 \sqrt{2}
   {\sigma_x}}\right)}
    \cos \left(\tfrac{\pi  m}{2}\right) \nonumber \\
& \times \exp\left(-\tfrac{m^2}{d_{yy}^2}\tfrac{\pi ^2}{L^2}  \tfrac {\sigma_x^2}{2 } \right) 
   \frac{\left(1  +(-1)^m \right)}{ m^2}    \exp\left( - \tfrac{m^2}{d^2_{yy}}  \tfrac{\pi^2}{L^2} D_{xx} t \right)  \Bigg]\\
M_{zz}(t) &= \frac{L^2}{12} \Bigg[1+
 \frac{24}{\pi^2}\sum_{n=1}^{n=\infty}  \frac{ 
 \Re\left[\text{erf}\left(\frac{{d^2_{zz}} L^2+2 i \pi  n \sigma_x^2}{2 \sqrt{2}
   L {{d_{zz}} \sigma_x}}\right)\right]}
   { \text{erf}\left(\frac{d_{zz} L}{2 \sqrt{2}
   {\sigma_x}}\right)}  
   \cos \left(\tfrac{\pi  n}{2}\right) \nonumber \\
&\times \exp\left(-\tfrac{n^2}{d_{zz}^2}\tfrac{\pi ^2}{L^2}  \tfrac {\sigma_x^2}{2 } \right)
    \frac{\left(1  +(-1)^n \right)}{ n^2}     \exp\left( -  \tfrac{n^2}{d^2_{zz}} \tfrac{\pi^2}{L^2} D_{xx} t \right)     \Bigg]
\end{align}
The normalized centered second moments are plotted for reference in Fig.~\ref{gaussian}.

\subsection{Case 4}

\begin{figure}
  \centering
 \includegraphics[scale=0.15]{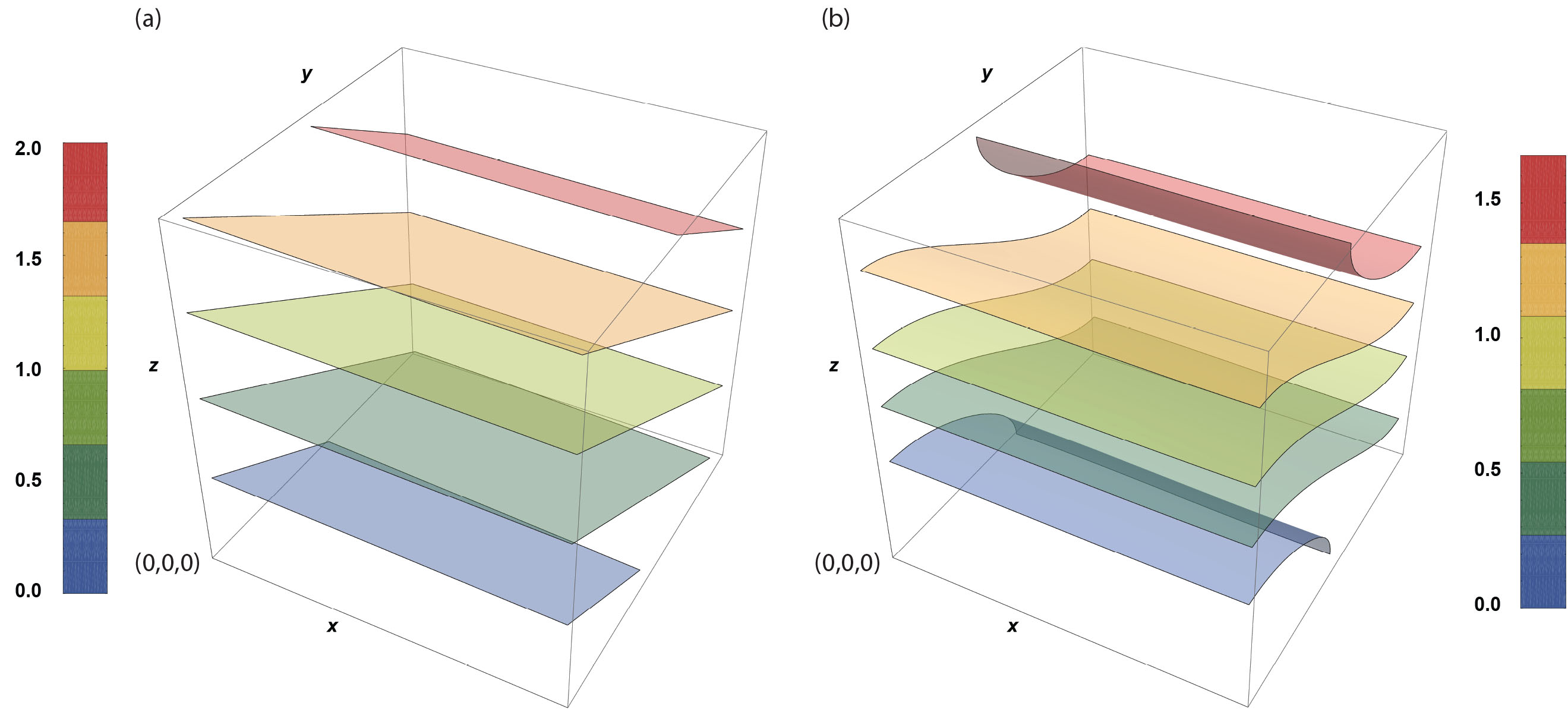}
  \caption{Concentration isosurfaces for the truncated Gaussian initial condition (Case 3);  concentration is normalized, $\tfrac{c}{c_\infty}$.  (a) The initial condition. (b) The concentration field at $t^*=\tfrac{1}{4} T^*_x$ using the first 20 terms of the series solution.  }
    \label{plane_concentration}
\end{figure}

The initial condition for Case 4 was selected specifically because it is not symmetric about the center of the domain.  Thus, the center of mass for the solution moves as a function of time.  Case 4 would be useful, for example, to provide a solution that would detect canceling errors in a code; symmetric initial conditions may miss such errors because of the symmetry imposed.   This case is also the only example provided in this paper for a non-multiplicative form for the initial condition.  Other non-multiplicative forms can be handled similarly.  To begin, we specify an initial condition that is a plane, given by
\begin{align}
&I.C.& c({\bf x},0)&= \frac{2}{L^4(1+\kappa_y+\kappa_z)}( x + \kappa_y y + \kappa_z  z)
\end{align}
Where $\kappa$ is a parameter that increases the slope of the plane along the $z-$ axis.  As with the previous initial conditions, for convenience the total mass for this function is normalized to unity.  Unlike previous cases, the evaluation of the Fourier coefficients are not simplified by separability.  Coefficients are computed by evaluating the integral given by Eq.~(\ref{fullFourier}).  Explicitly, this is

\begin{equation}
 \bar{B}_{\ell m n} =  \frac{2}{L^4(2+a)} \frac{2^{N_{0}}}{L^3}  \int_{x=0}^{x=L} \int_{y=0}^{y=L} \int_{z=0}^{z=L} ( x + \kappa_y y + \kappa_z z) \cos(\ell \pi \tfrac{x}{L})  \cos(m \pi \tfrac{y}{L}) \cos(n \pi \tfrac{z}{L}) \;dx \;dy  \;dz
\end{equation}
where $\ell, m, n = 0, 1, 2, \ldots $.  A simple integration shows that $\bar{B}_{\ell m n}$ is zero when any two of the indexes are greater than zero.  Thus, the only non-zero contributions for the coefficients are of the form $B_{\ell 0 0}$, $B_{0 m 0}$, and $B_{0 0 n}$.  This yields
 
 \begin{equation}
 \bar{B}_{\ell 0 0} = \frac{2^{N_0+1}}{L^5(1+\kappa_y+\kappa_z)}  \int_{x=0}^{x=L} x \cos(\ell \pi \tfrac{x}{L})\,dx = \left\{ {\begin{array}{*{20}{c}}
  \frac{1}{(1+\kappa_y+\kappa_z)}{\frac{4(-1+(-1)^{\ell})}{\ell^2  \pi^2}\frac{1}{L^3},\,\,\, for~\ell>0} \\ 
  {\frac{1}{(1+\kappa_y+\kappa_z)}\frac{1}{L^3}, \,\,\, for~\ell=0 }
\end{array}} 
 \right.
 \end{equation}
 \begin{equation}
 \bar{B}_{0 m 0} = \frac{2^{N_0+1}}{L^5(1+\kappa_y+\kappa_z)}  \int_{y=0}^{y=L} \kappa_y y \cos(m \pi \tfrac{y}{L})\,dy = \left\{ {\begin{array}{*{20}{c}}
  \frac{\kappa_y}{(1+\kappa_y+\kappa_z)} {\frac{4(-1+(-1)^{m})}{m^2  \pi^2}\frac{1}{L^3},\,\,\, for~m>0} \\ 
  {\frac{\kappa_y}{(1+\kappa_y+\kappa_z)}\frac{1}{L^3}, \,\,\, for~m=0 }
\end{array}} 
 \right.
 \end{equation}
 \begin{equation}
 \bar{B}_{0 0 n} = \frac{2^{N_0+1}}{L^5(1+\kappa_y+\kappa_z)}  \int_{z=0}^{z=L} \kappa_z z \cos(m \pi \tfrac{z}{L})\,dz = \left\{ {\begin{array}{*{20}{c}}
  \frac{\kappa_z}{(1+\kappa_y+\kappa_z)}{\frac{4\kappa(-1+(-1)^{n})}{n^2 \pi^2}\frac{1}{ L^3},\,\,\, for~n>0} \\ 
  \frac{\kappa_z}{(1+\kappa_y+\kappa_z)}{\frac{1}{L^3}, \,\,\, for~n=0 }
\end{array}} 
 \right.
 \end{equation}
 Recall, the general solution for the non-multiplicative case is given by Eq.~(\ref{generalSol}); upon substitution, this yields
 
 \begin{align}
c_{A\gamma}(x,y,z,t) =& \frac{1}{L^3}\Bigg[1+
  \frac{1}{(1+\kappa_y+\kappa_z)}\sum_{\ell = 1}^{\ell=\infty}  \exp \left( -  \ell^2  \tfrac{\pi^2}{L^2} D_{xx}  t \right) 
 \frac{4(-1+(-1)^{\ell})}{\ell^2  \pi^2} \cos(\ell \pi \tfrac{x}{L}) \nonumber \\
 &+\frac{\kappa_y}{(1+\kappa_y+\kappa_z)} \sum_{m = 1}^{m=\infty}  \exp \left( -  \tfrac{m^2}{d^2_{yy}}  \tfrac{\pi^2}{L^2} D_{xx}  t \right) 
 \frac{4(-1+(-1)^{m})}{m^2  \pi^2} \cos(m \pi \tfrac{y}{L}) \nonumber \\
  &+\frac{\kappa_z}{(1+\kappa_y+\kappa_z)}\sum_{n = 1}^{n=\infty}    \exp \left( -\tfrac{n^2}{d^2_{zz}}  \tfrac{\pi^2}{L^2} D_{xx}   t \right) 
\frac{4(-1+(-1)^{n})}{n^2  \pi^2} \cos(n \pi \tfrac{y}{L})     
 \Bigg] 
 \label{weird}    
\end{align}
Note that here, we combined the constant term outside the sum to be consistent with the final form of the previous solutions.  As a comment, also note that this problem is the superposition of the three analogous problems in 1-dimension.  This is to be expected because of the linearity of the diffusion equation, the boundary conditions, and the initial condition.  If the initial condition had not been linear, superposition would not have been possible, and the result would have been much more complex.  Isosurface plots of the normalized concentration appear in Fig.~\ref{plane_concentration}.

The moments for this initial condition are given by the following.\\

{\noindent}{\it Zeroth moment}
\begin{equation}
m_0(t)  = \int_{x=0}^{x=L}  \int_{y=0}^{y=L}  \int_{z=0}^{z=L}  c_{A\gamma}(x,y,z,t) dz dy dx = 1
\end{equation}

 {\noindent}{\it First moments}
 Unlike the previous cases, the first moments for this initial condition are not constant.  Integrating Eq.~(\ref{weird}) directly (via Eqs.~(\ref{firstx})-(\ref{firstz})), the first moments are given by
 
 \begin{align}
m_x(t) =& \frac{L}{2} \Bigg[ 1+
\frac{1}{(1+\kappa_y+\kappa_z)} \sum_{\ell = 1}^{\ell=\infty}  \exp \left( -  \ell^2  \tfrac{\pi^2}{L^2} D_{xx}  t \right) 
 \frac{8(-1+(-1)^{\ell})^2}{\ell^4  \pi^4} 
 \Bigg]  \label{firstmo1}  \\
 m_y(t)=& \frac{L}{2}\Bigg[ 1+
\frac{\kappa_y}{(1+\kappa_y+\kappa_z)} \sum_{m = 1}^{m=\infty}  \exp \left( -  \tfrac{m^2}{d^2_{yy}}  \tfrac{\pi^2}{L^2} D_{xx}  t \right) 
 \frac{8(-1+(-1)^{m})^2}{m^4  \pi^4} 
 \Bigg] \\
 m_z(t)= & \frac{L}{2}\Bigg[1+
\frac{ \kappa_z}{(1+\kappa_y+\kappa_z)} \sum_{n = 1}^{n=\infty}  \exp \left( -\tfrac{n^2}{d^2_{zz}}  \tfrac{\pi^2}{L^2} D_{xx}   t \right) 
 \frac{8(-1+(-1)^{n})^2}{n^4 \pi^4}  \Bigg]  
 \label{firstmo3}
\end{align}

\begin{figure}
  \centering
 \includegraphics[scale=0.8]{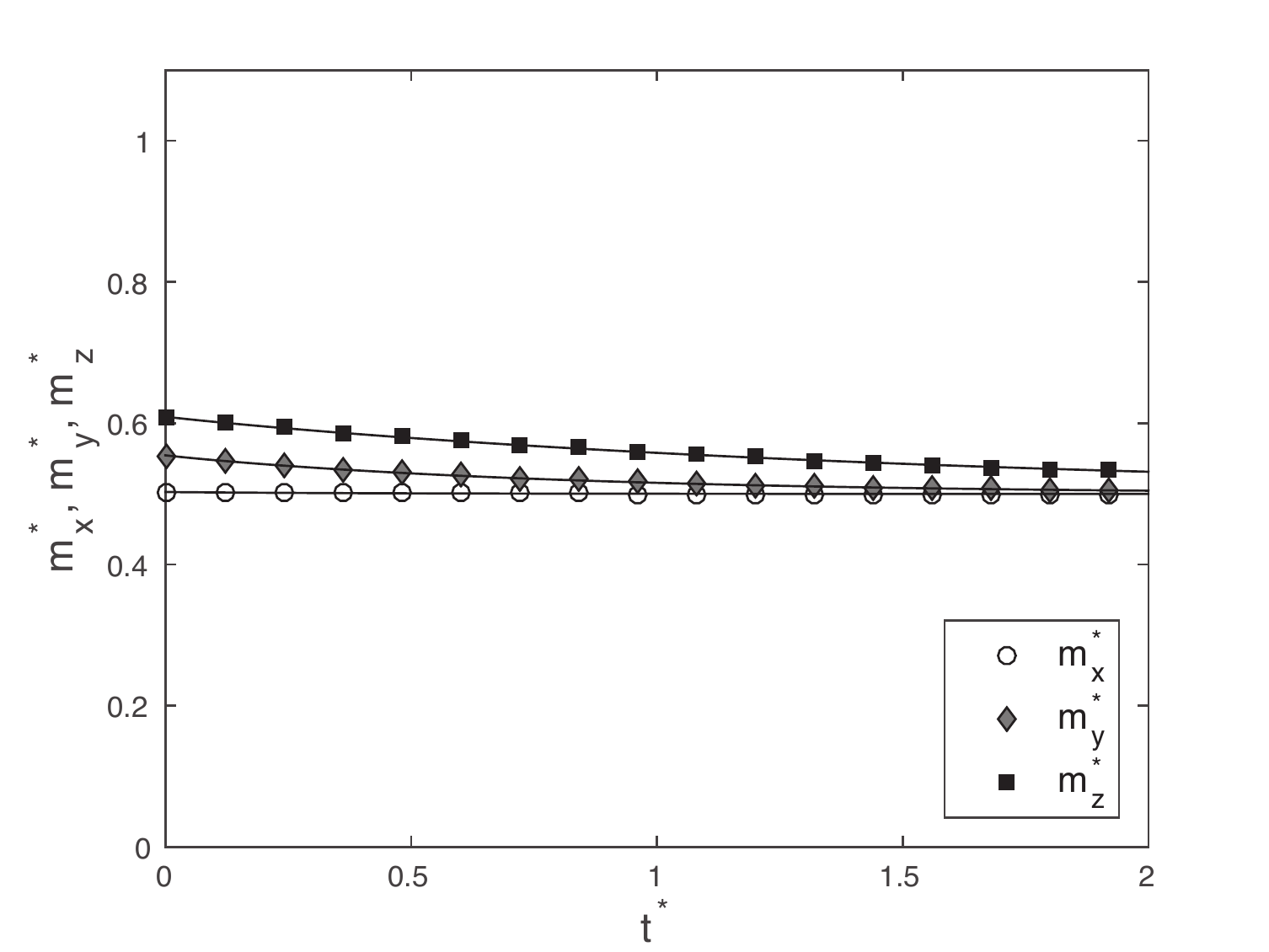}
  \caption{First moments for the plane initial condition (Case 4) via series solution (solid lines), and numerical computations  (symbols).}
      \label{plane1st}
\end{figure}

The first moments in this case are different from the first three cases because the center of mass of the initial condition is different from the center of mass of the steady-state condition; thus, the first moments evolves in time.  A plot of the first moments as a function of time is provided in Fig.~\ref{plane1st}.  Note for this plot, we have defined the dimensionless variables

\begin{align}
m^*_{x}=& \frac{m_{x}}{L }   \\
m^*_{y}=& \frac{m_{y}}{L }  \\
m^*_{z}=& \frac{m_{z}}{L }  
\end{align}

For this case, the centered second moments are complicated by the transient first moments

\begin{align}
M_{xx}(t) &=  \int_{x=0}^{x=L}  \int_{y=0}^{y=L}  \int_{z=0}^{z=L}    (x-m_x(t))^2 c_{A\gamma}(x,y,z,t) dz dy dx \label{secmo1}  \\
M_{yy}(t) &=  \int_{x=0}^{x=L}  \int_{y=0}^{y=L}  \int_{z=0}^{z=L}    (y--m_y(t))^2 c_{A\gamma}(x,y,z,t) dz dy dx   \\
M_{zz}(t) &=  \int_{x=0}^{x=L}  \int_{y=0}^{y=L}  \int_{z=0}^{z=L}    (z--m_z(t))^2 c_{A\gamma}(x,y,z,t) dz dy dx  \label{secmo3}
\end{align}
Substituting Eqs.~(\ref{weird}) and (\ref{firstmo1})-(\ref{firstmo1}) into Eqs.~(\ref{secmo1})-(\ref{secmo1}), the closed-form expressions for the centered second moment are \\

\begin{figure}
  \centering
 \includegraphics[scale=0.8]{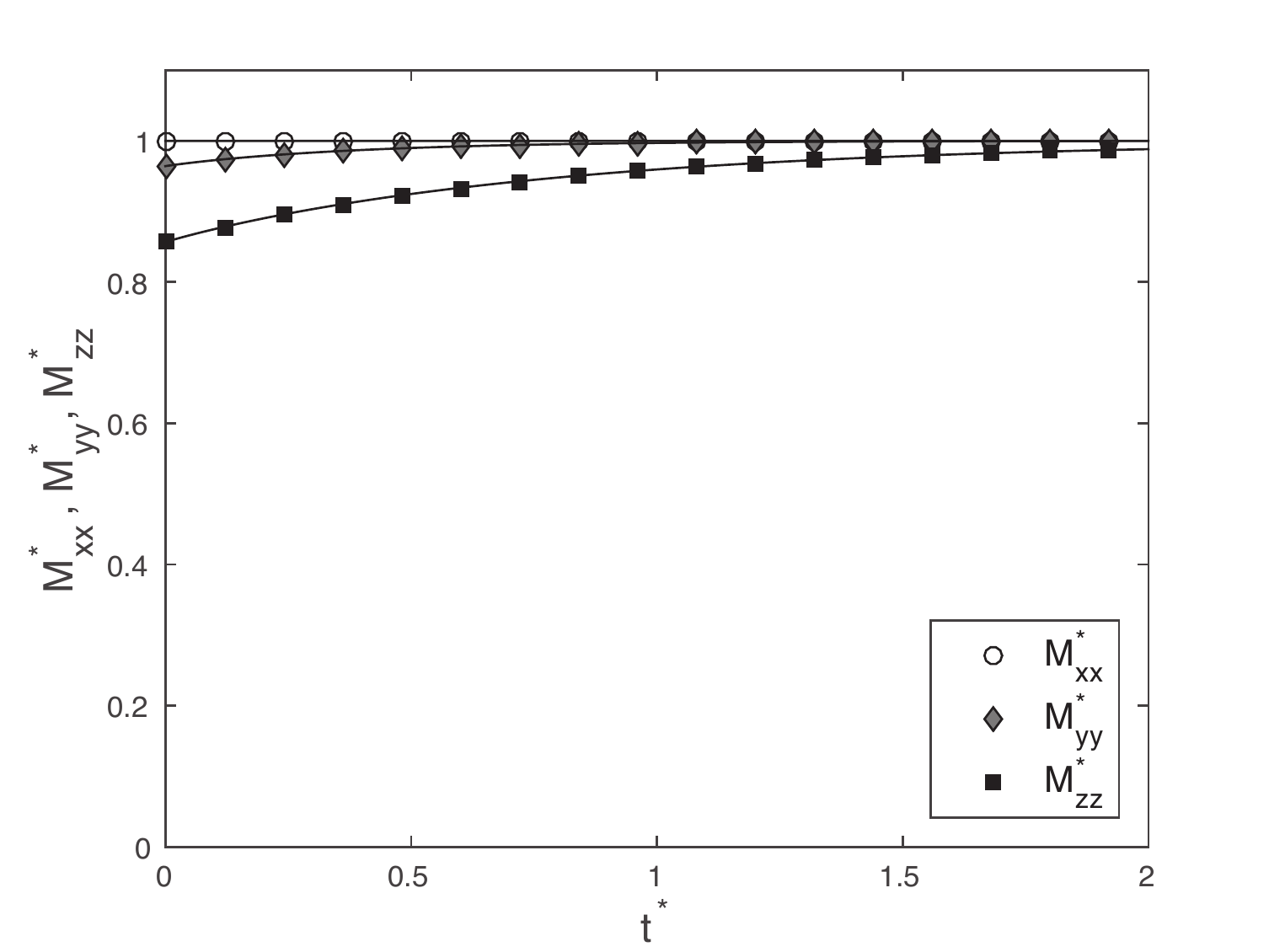}
  \caption{Centered second moments for the plane initial condition (Case 4) via series solution (solid lines), and numerical computations  (symbols).}
      \label{plane}
\end{figure}

{\noindent}{\it Centered second moments}
\begin{align}
M_{xx}(t) &=  \frac{L^2}{3}-L m_x(t) +m^2_x(t) \nonumber \\
&+ \frac{1}{(1+\kappa_y+\kappa_z)} L^2  \sum_{\ell = 1}^{\ell=\infty} \exp\left( -\ell^2 \tfrac{\pi^2}{L^2} D_{xx} t \right) \frac{8(-1+(-1)^\ell)}{\ell^4 \pi^4}{(-1)^\ell} \nonumber \\
&-\frac{1}{(1+\kappa_y+\kappa_z)} L m_x(t)  \sum_{\ell = 1}^{\ell=\infty} \exp\left( -\ell^2 \tfrac{\pi^2}{L^2} D_{xx} t \right) \frac{8(-1+(-1)^\ell)^2}{\ell^4 \pi^4}
\label{mx99}
\end{align}
\begin{align}
M_{yy}(t) &=  \frac{L^2}{3}-L m_y(t) +m^2_y(t) \nonumber \\
&+  \frac{\kappa_y}{(1+\kappa_y+\kappa_z)} L^2  \sum_{m = 1}^{m=\infty} \exp\left( -\tfrac{m^2}{d^2_{yy} } \tfrac{\pi^2}{L^2} D_{xx} t \right) \frac{8(-1+(-1)^m)}{m^4 \pi^4}{(-1)^m} \nonumber \\
&- \frac{\kappa_y}{(1+\kappa_y+\kappa_z)} L m_y(t)  \sum_{m = 1}^{m=\infty} \exp\left( -\tfrac{m^2}{d^2_{yy} } \tfrac{\pi^2}{L^2} D_{xx} t \right) \frac{8(-1+(-1)^\ell)^2}{\ell^4 \pi^4}
\label{my99}
\end{align}
\begin{align}
M_{zz}(t) &=  \frac{L^2}{3}-L m_z(t) +m^2_z(t) \nonumber \\
&+  \frac{\kappa_z}{(1+\kappa_y+\kappa_z)} L^2  \sum_{n = 1}^{n=\infty} \exp\left(-\tfrac{n^2}{d^2_{zz} } \tfrac{\pi^2}{L^2} D_{xx} t \right) \frac{8(-1+(-1)^n)}{n^4 \pi^4}{(-1)^n} \nonumber \\
&- \frac{\kappa_z}{(1+\kappa_y+\kappa_z)} L m_x(t)  \sum_{n = 1}^{n=\infty} \exp\left( -\tfrac{n^2}{d^2_{zz} } \tfrac{\pi^2}{L^2} D_{xx} t \right) \frac{8(-1+(-1)^n)^2}{n^4 \pi^4} 
\label{mz99}
\end{align}
As a check on this result, one can observe that as time becomes arbitrarily large, each of $m_x$, $m_y$, and $m_z$ tend toward the value $L/2$.  Thus, each of the second centered moments tend the to the value $\tfrac{L^2}{3}-\tfrac{L^2}{2}+\tfrac{L^2}{4}= \tfrac{L^2}{12}$, which the correct value for the steady-state centered second moment.  A plot of the centered second moments appear in Fig.~\ref{plane}.

\section{Computations} \label{computations}

To demonstrate the utility of the analytical solutions, we conducted direct numerical simulations of diffusion in a cube for a prescribed anisotropic effective diffusivity tensor under the four different initial conditions using the finite element method (FEM) as coded in the
commercial package COMSOL Multiphysics 4.4. This code has been previously verified to
be a correct FEM representation of the diffusion/heat equations (e.g., \citep{xu2010verification}), in the sense of
\citet{roache1998verification}. For the simulations in this paper, the parameters summarized in Table \ref{table-parameters} were
adopted. The code COMSOL solves the discretized equations using the generalized minimal
residual method with geometric multigrid preconditioning. Successive over-relaxation was
used in pre- and post-smoothing. To assure that the numerical results were converged, the
domain was tessellated using a tetrahedral mesh at three levels of refinement.\\

We use the Grid Convergence Index (GCI) as proposed by \citet{roache_1994} to estimate the discretization uncertainties.  
The GCI is very closely based on Richardson extrapolation and has been developed to serve as an uncertainty and convergence analysis
tool in computational fluid dynamics when analytical solutions are not available. The basic idea is to 
obtain the solution for a given variable $\phi$ on (at least) three grids with spacings $\Delta_{1}$, $\Delta_{2}$, and $\Delta_{3}$, for which
the refinement ratios $r_{21} = \Delta_{2}$ / $\Delta_{1}$, and $r_{32} = \Delta_{3}$ / $\Delta_{2}$, and the variations 
$\epsilon_{32}(\textbf{x}) = \phi_{3}(\textbf{x}) - \phi_{2}(\textbf{x})$, and $\epsilon_{21}(\textbf{x}) = \phi_{2}(\textbf{x}) - \phi_{1}(\textbf{x})$ are 
calculated for values of the fine and medium grid solutions ($\phi_{1}$ and $\phi_{2}$) interpolated on the coarse grid. A ${\it local}$ (i.e. pointwise) 
order of accuracy $\mathscr{P}(\textbf{x})$ can be obtained by solving the following equations iteratively \cite{celik_2008}:

\begin{equation}
 \mathscr{P}(\textbf{x}) = \frac{1}{\mathit{ln}(r_{21})} \left \lvert \mathit{ln} \left \lvert \frac{\epsilon_{32}(\textbf{x})}{\epsilon_{21}(\textbf{x})} \right \rvert + \mathit{ln} \left ( \frac{r_{21}^{\mathscr{P}(\textbf{x})} - s(\textbf{x})}{r_{32}^{\mathscr{P}(\textbf{x})} - s(\textbf{x})} \right) \right \rvert
\end{equation}
 
\begin{equation}
 s(\textbf{x}) = \begin{cases}
                  +1, \quad \mbox{if} \quad \epsilon_{32}(\textbf{x}) \; \epsilon_{21}(\textbf{x}) > 0 \\
		  -1, \quad \mbox{if} \quad \epsilon_{32}(\textbf{x}) \; \epsilon_{21}(\textbf{x}) < 0 \\
                 \end{cases}
\end{equation}

The fine-grid convergence index is then computed using:

\begin{equation}
  \label{gci}
 \mbox{GCI}_{\mathit{fine}}^{21} = \frac{1.25 \; e_{a}^{21}}{r_{21}^{\mathscr{P}} - 1},
\end{equation}

where $e_{a}^{21} = \lvert (\phi_{1} - \phi_{2})$ / $\phi_{1} \rvert$. \\

\begin{figure}
   \centering
 \includegraphics[scale=0.8]{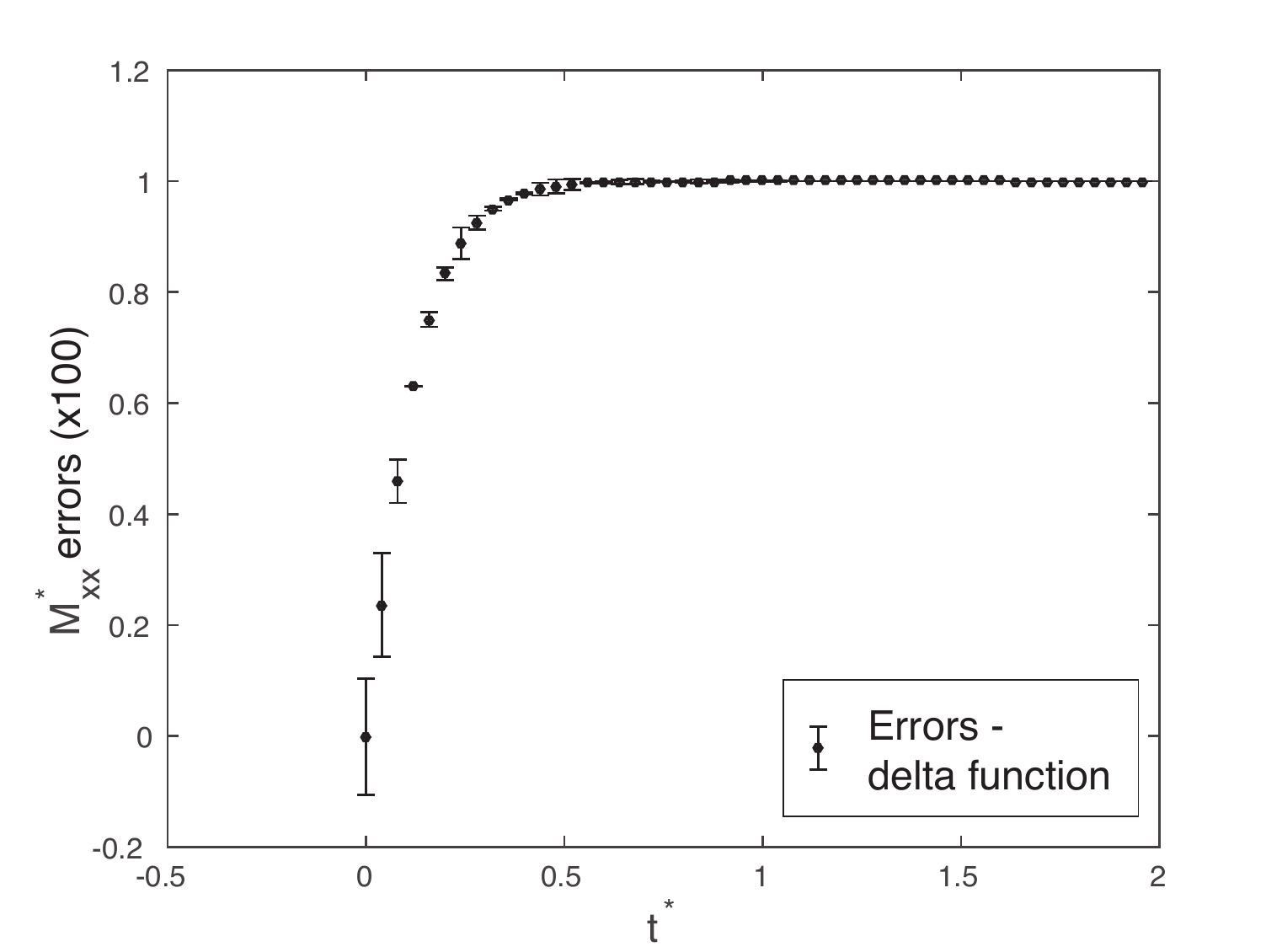}
  \caption{Sample discretization uncertainty estimates for $M_{xx}$ for the $\delta$-function initial condition; the other three initial conditions exhibited similar local errors. Error bars have been magnified 100x to facilitate inspection.}
   \label{error}
\end{figure}

Equation (\ref{gci}) can be used to estimate the discretization uncertainty locally or globally. $\mathscr{P}$ and $\mbox{GCI}_{\mathit{fine}}^{21}$ acquire single 
values for a global variable (e.g. total mass). When local uncertainty estimates are desired, the GCI can be 
computed locally using a global average of $\mathscr{P}$ \cite{celik_2008} as demonstrated in Fig.~\ref{error}. Based on these, 
the numerical uncertainties in the fine-grid solutions for the total mass, $m_{0}$, were 0.62$\%$, 0.56$\%$, and 2.06$\%$ for the 
$\delta$-function, Gaussian, and step function cases, respectively (data not shown).  We define the percentage local discretization uncertainty for the second moment by

\begin{equation}
\mu_2 = \frac{GCI^{21}_{fine}}{M_\infty}\times 100
\end{equation}
where in this case $GCI^{21}_{fine}$ is the value obtained for computed time point for applicable second moment, $M_{xx}, M_{yy}$, or, $M_{zz}$. For all for cases the maximum value of $\mu_2$ for all three moments was less than 10$\%$.

\section*{Acknowledgements}
This work was supported in part by the National Science Foundation, under grant CBET-336983.



%

%
\appendix
\renewcommand{\thesection}{\Alph{section}}
\renewcommand{\theequation}{\Alph{section}.\arabic{equation}}
\renewcommand\thefigure{\thesection.\arabic{figure}} 
\renewcommand\thetable{\thesection.\arabic{table}} 
\setcounter{figure}{0}
\setcounter{equation}{0}
\setcounter{table}{0}

\section{Appendix}\label{appendix}

In this appendix, we briefly explain the transformation of variables required for orthotropic heat or diffusion problems in a rectangular parallelpiped.  Coordinate vectors in this system are specified by ${\bar{\bf x}}=(\bar{x}, \bar{y}, \bar{z})$, and the parallelpiped has sides of length $L$ in the $\bar{x}-$direction, and $L_{\bar{y}}$, and $L_{\bar{z}}$ the the other two perpendicular directions.   For this case, the problem is specified by

\begin{align}
&&\frac{\partial c_{A\gamma}}{\partial t} =D_{\bar{x}\bar{x}} \frac{\partial^2 c_{A\gamma}}{\partial \bar{x}^2} + & D_{\bar{y}\bar{y}} \frac{\partial^2 c_{A\gamma}}{\partial \bar{y}^2} +D_{\bar{z}\bar{z}} \frac{\partial^2 c_{A\gamma}}{\partial \bar{z}^2} \\
&B.C.~1& -{\bf n}_{\gamma} \cdot  {\textbf{\sffamily \bfseries D}}_{\gamma}\cdot \nabla c_{A\gamma}  &= 0,~~\left\{ \begin{gathered}
  at~all~boundaries \\ 
  of~the~rectangular \\ 
  parallelpiped \\ 
\end{gathered}  \right. \\
&I.C.& c(\bar{{\bf x}},0)&=\bar{\Phi}(\bar{x}, \bar{y}, \bar{z})
\end{align}
Now, define the new variables

\begin{align}
x & = \bar{x}   \\
y & = {\bar y}\frac{L_{\bar{y}}}{L}   \\
z & = {\bar z}\frac{L_{\bar{z}}}{L}
\end{align}
From these definitions, it is clear that coordinate vectors $(x,y,z)$ span $(0,0,0)$ to $(L, L, L)$ as the coordinate vectors $(\bar{x},\bar{y},\bar{z})$ span $(0,0,0)$ to $(L, L_{\bar{y}}, L_{\bar{z}})$. 

Applying the chain rule twice, it is easy to determine the relationships

\begin{align}
\frac{\partial^2 c_{A\gamma}} {\partial^2 \bar{x}} & = \frac{\partial^2 c_{A\gamma}}{\partial x^2}   \\
\frac{\partial^2 c_{A\gamma}} {\partial^2 \bar{y}} & = \frac{L^2} {L^2_{\bar{y}}}     \frac{\partial^2 c_{A\gamma}}{\partial y^2}   \\
\frac{\partial^2 c_{A\gamma}} {\partial^2 \bar{z}} & = \frac{L^2}{L^2_{\bar{z}}}     \frac{\partial^2 c_{A\gamma}}{\partial z^2}   
\end{align}
Finally, denoting

\begin{align}
D_{xx} = D_{\bar{x}\bar{x}} \\
D_{yy} = D_{\bar{y}\bar{y}} \frac{L^2} {L^2_{\bar{y}}}  \\
D_{zz} = D_{\bar{z}\bar{z}} \frac{L^2} {L^2_{\bar{z}}}  \\
\end{align}
we see that the the final form of the governing equation, boundary conditions, and initial conditions takes the form

\begin{align}
&&\frac{\partial c_{A\gamma}}{\partial t} =D_{{x}{x}} \frac{\partial^2 c_{A\gamma}}{\partial {x}^2} + & D_{{y}{y}} \frac{\partial^2 c_{A\gamma}}{\partial {y}^2} +D_{{z}{z}} \frac{\partial^2 c_{A\gamma}}{\partial {z}^2} \\
&B.C.~1& -{\bf n}_{\gamma} \cdot  {\textbf{\sffamily \bfseries D}}_{\gamma}\cdot \nabla c_{A\gamma}  &= 0,~~\begin{array}{*{20}{c}}
  {at~all~boundaries} \\ 
  {of~the~cube} \\
\end{array}   \\
&I.C.& c({\bf x},0)&=\Phi(x,y,z)
\end{align}
where 
$\Phi(x,y,z) = \bar{\Phi}({x}, \tfrac{L} {L_{\bar{y}} }y , \tfrac{L}{L_{\bar{z}} } z )$
represents the initial condition in the new coordinate system.  This problem is identical to that in the main body of the paper, establishing that for all problems on rectangular parallelpiped, there is an equivalent problems for a cubic domain.

\bibliographystyle{elsarticle-num-names}
\bibliography{diffuse}

%
%
%
\end{document}